\documentclass[11pt]{article}
\usepackage{amsfonts}
\usepackage{latexsym}

\newcommand{\qed}{$\hfill\Box$}

\newcommand{\R}{\mathbb{R}}
\newcommand{\LL}{\mathbb{L}}
\newcommand{\E}{\mathbb{E}}
\newcommand{\PP}{\mathbb{P}}
\newcommand{\N}{\mathbb{N}}

\newcommand{\al}{\alpha}
\newcommand{\be}{\beta}
\newcommand{\la}{\lambda}
\newcommand{\ga}{\gamma}
\newcommand{\ka}{\kappa}
\newcommand{\Ga}{\Gamma}
\newcommand{\si}{\sigma}
\newcommand{\Si}{\Sigma}
\newcommand{\ep}{\varepsilon}
\newcommand{\vpi}{\varpi}

\newcommand{\de}{\delta}
\newcommand{\te}{\theta}

\newcommand{\De}{\Delta}
\newcommand{\om}{\omega}
\newcommand{\Om}{\Omega}
\newcommand{\ze}{\zeta}
\newcommand{\vp}{\varphi}

\newcommand{\n}{\mathcal N}

\newcommand{\ea}{\mathcal E}
\newcommand{\g}{\mathcal G}

\newcommand{\f}{\mathcal F}
\newcommand{\rR}{\mathcal R}

\newcommand{\vv}{\mathcal V}
\newcommand{\lll}{\mathcal L}

\def\nib{\noindent\bf}
\def\ni{\noindent}

\newcommand{\rn}{\sqrt{\De_n}}

\newcommand{\st}{\sum_{i=1}^{[t/\De_n]}}
\newcommand{\ec}{\mathbb{E}^n_{i-1}}
\newcommand{\pc}{\mathbb{P}^n_{i-1}}
\newcommand{\dd}{\Delta^n_i}
\newcommand{\is}{\star (\mu -\nu)}
\newcommand{\itai}{\int^{i\De_n}_{(i-1)\De_n}}
\newcommand{\proba}{(\Omega ,\f,(\f_t)_{t\geq0},\PP)}
\def\probt{\hbox{$(\widetilde{\Omega},\widetilde{\f},
(\widetilde{\f}_t)_{t\geq0},\widetilde{\PP})$}}

\newcommand{\toop}{\stackrel{\PP}{\longrightarrow}}
\newcommand{\tols}{~\stackrel{\lll-(s)}{\longrightarrow}~}

\newcommand{\toSp}{\stackrel{\mbox{\tiny Sk.p.}}{\longrightarrow}}
\newcommand{\tosc}{\stackrel{\mbox{\tiny Sk}}{\longrightarrow}}
\newcommand{\toucp}{~\stackrel{\mbox{\tiny u.c.p.}}{\longrightarrow}}
\newcommand{\tovp}{\stackrel{\mbox{\tiny v.p.}}{\longrightarrow}}
\newcommand{\tOu}{\mbox{\rm O}_{Pu}}
\newcommand{\tou}{\mbox{\rm o}_{Pu}}

\newcommand{\wbe}{\widehat{\beta}}
\newcommand{\wxi}{\widehat{\xi}}

\newcommand{\wga}{\widehat{\ga}}

\newcommand{\WOm}{\widetilde{\Omega}}

\newcommand{\WP}{\widetilde{\mathbb{P}}}
\newcommand{\WE}{\widetilde{\mathbb{E}}}
\newcommand{\WQ}{\widetilde{Q}}
\newcommand{\Wf}{\widetilde{\f}}

\newcommand{\Wsi}{\widetilde{\sigma}}
\newcommand{\Wga}{\widetilde{\gamma}}
\newcommand{\Wde}{\widetilde{\delta}}
\newcommand{\Wxi}{\widetilde{\xi}}
\newcommand{\Wb}{\widetilde{b}}

\newcommand{\BH}{\overline{H}}

\newcommand{\BB}{\overline{B}}
\newcommand{\BU}{\overline{U}}
\newcommand{\BV}{\overline{V}}
\newcommand{\BW}{\overline{W}}

\newcommand{\BK}{\overline{K}}
\newcommand{\Bg}{\overline{g}}
\newcommand{\BM}{\overline{M}}
\newcommand{\BX}{\overline{X}}
\newcommand{\Bk}{\overline{k}}

\newcommand{\umu}{\underline{\mu}}
\newcommand{\unu}{\underline{\nu}}
\newcommand{\uF}{\underline{F}}

\newcommand{\vsc}{\vskip 5mm}

\newcommand{\vst}{\vskip 3mm}
\newcommand{\vsq}{\vskip 4mm}

\newcommand{\bee}{\begin{equation}}
\newcommand{\eee}{\end{equation}}
\newcommand{\bea}{\begin{eqnarray}}
\newcommand{\eea}{\end{eqnarray}}
\newcommand{\bean}{\begin{eqnarray*}}
\newcommand{\eean}{\end{eqnarray*}}

\renewcommand{\theequation}{\arabic{section}.\arabic{equation}}

\newtheorem{prop}{Proposition}[section]

\newtheorem{lem}[prop]{Lemma}
\newtheorem{theo}[prop]{Theorem}
\newtheorem{rem}[prop]{Remark}

\begin{document}

\title{Asymptotic properties of realized power variations and
  related functionals of semimartingales}
\author{Jean Jacod
\thanks{Institut de Math\'ematiques de Jussieu, 175 rue du Chevaleret
75 013 Paris, France (CNRS -- UMR 7586,
and Universit\'e Pierre et Marie Curie - P6)}}

%\date{}

\maketitle

\nib Abstract. \rm\small This paper is concerned with the asymptotic
behavior of sums of the form $U^n(f)_t=\st f(X_{i\De_n}-X_{(i-1)\De_n})$, where
$X$ is a $1$-dimensional semimartingale and $f$ a suitable test
function, typically $f(x)=|x|^r$, as $\De_n\to0$. We prove a variety
of ``laws of large 
numbers'', that is convergence in probability of $U^n(f)_t$, sometimes
after normalization. We also exhibit in many cases the rate of
convergence, as well as associated central limit theorems. \normalsize

\noindent AMS classification\,: 60F17, 60G48

\noindent Keywords\,: Central limit theorem, quadratic variation,
power variation, semimartingale.

\normalsize

\section{Introduction}\label{sec-Intro}
\setcounter{equation}{0}
\renewcommand{\theequation}{\thesection.\arabic{equation}}

In many practical situations one observes a process $X$ at finitely
many times, and from these observations one wants to infer various properties
of the process. For example, in finance the price of an
asset is observed at discrete times and one aims to determine the
volatility or the integrated volatility, or perhaps the presence of
jumps and some properties about their sizes. In statistics one wants
to determine the parameters on which the law of the process depends,
or one may want to perform some non-parametric inference on the model.

There are indeed two very different situations. One is when the
observations occur at time $~0,\De,2\De,\cdots,n\De~$ for a fixed time
lag $\De$, whereas $n$ is large\,: then any kind of inference
necessitates some ``ergodic'' properties of the basic process. Another
situation is what is called {\em high frequency}\, observations, where the
time lag $\De$ is small, which in the asymptotic setting means that we
let $\De=\De_n$ depend on the number $n$ of observations and go to $0$
as $n\to\infty$. This second situation is the one we are interested
in here.

A first, well known, example of how discrete observations of $X$
allow to approximate some basic characteristics of the process is the
convergence of the ``realized'' (or approximate) quadratic variation
towards the ``true'' one. More generally one may look at the realized
$r$-th power variation at stage $n$, that is the (observable) process
\bee\label{I0}
\{X\}^{r,n}_t=\st|X_{i\De_n}-X_{(i-1)\De_n}|^r.
\eee
When $r=2$ the processes $\{X\}^{2,n}$ converge (as $\De_n\to0$)
to $[X,X]$, the quadratic variation of $X$, as soon as $X$ is a
semimartingale, and even in some more general situations. When $r>2$
then $\{X\}^{r,n}_t$ converges to $\sum_{s\leq t}|\De X_s|^r$ (where
$\De X_s$ is the size of the jump of $X$ at time $s$) for any
semimartingale, also an old result due to L\'epingle in
\cite{L}. When $r\in(0,2)$ then $\{X\}^{r,n}_t$ blows up in general,
but $\De_n^{1-r/2}\{X\}^{r,n}_t$ converges to the continuous part of
$[X,X]_t$\,: this does not hold in general, though, but under some
(weak) assumptions on $X$. So this allows in principle to ``separate''
the jumps of $X$ from its continuous part.

Again for practical applications, having the convergence of
$\{X\}^{r,n}$ (possibly after normalization) is not enough, we need
rates and, if possible, an associated central limit theorem. This
describes the main aim of this paper: find conditions for the above
convergence, and for associated CLTs when they exist. We do that for the
processes $\{X\}^{r,n}$, and more generally for the following
processes, for suitable test functions $f$ and cut-off exponent
$\vpi>0$ and level $\al>0$\,:
\bee\label{I4}
\left.\begin{array}{l}
V^n(f)_t= \st f(X_{i\De_n}-X_{(i-1)\De_n}), \\[2mm]
V'^n(f)_t=\st f((X_{i\De_n}-X_{(i-1)\De_n})/\rn),\\[2mm]
V''^n(\vpi,\al)_t=\st(X_{i\De_n}-X_{(i-1)\De_n})^2
1_{\{|X_{i\De_n}-X_{(i-1)\De_n}|\leq\al\De_n^{\vpi}\}}\,.
\end{array}\right\}
\eee
The convergence of these processes and the associated CLTs
hold or not, depending on the properties of $f$ of 
course, and especially on its behavior near $0$, but also on the
properties of the basic semimartingale $X$. Note that we {\it always
  assume that}\, $\De_n\to0$.

The reader may find motivations and practical uses of realized power
variations in finance in Andersen, Bollersley and Diebold \cite{ABD}
or Barndorff-Nielsen and Shephard \cite{BS} and references therein,
for continuous processes. The later authors have also introduced and
thoroughly used the ``bi-power 
variations'' where the summands in (\ref{I0}) are products of powers
of two successive increments of $X$ instead of one, and probably what
follows can also be done for bi- or multi-power variations as
well. The case where $X$ is discontinuous has been studied by Mancini
\cite{Ma1}, \cite{Ma2} (using processes similar to $V''^n(\vpi,\al)$)
and Woerner 
\cite{W1}, \cite{W2} (for the power and bi-power variations) and recently by
Barndorff-Nielsen, Shephard and Winkel \cite{BSW}, and in those papers
special cases of the forthcoming results may be found. 

In \cite{J2} we have considered the same problems than here when $X$
is a L\'evy processes, with almost complete answers. 
In the semimartingale case the picture shown below
is neither as good nor as complete as in the L\'evy case. The
proofs are mostly quite different (except for Theorem \ref{T1}), hence
this paper is essentially
independent of \cite{J2} although the basic ideas are the same. On the
other hand, some of the results here heavily
rely upon the paper \cite{BGJPS} in which similar problems
have been solved when $X$ is {\em continuous}. 

Let us also mention that only the $1$-dimensional case is
considered here, although it covers the case where $X$ is one of the
components of a multidimensional semimartingale. Some results
obviously hold as well when $X$ is multidimensional (those concerned
with $V'^n(f)$ in particular), others do not: if 
$f$ is singular at $0$, the description of the singularity in the
multidimensional case is clearly much more sophisticated than in
dimension $1$.

The main notation, assumptions and results are gathered in Section
\ref{secR}. All (unfortunately rather tedious) proofs are in
the subsequent sections.

\section{Notation, assumptions, results}\label{secR}
\setcounter{equation}{0}
\renewcommand{\theequation}{\thesection.\arabic{equation}}

\subsection{Some general notation.}

Let us first introduce a number of notation to be used throughout.
With any process $Y$ we associate its increments $\dd Y$ and the
''discretized process'' as follows
\bee\label{I2}
\dd Y=Y_{i\De_n}-Y_{(i-1)\De_n},\qquad
Y^{(n)}_t~=~Y_{\De_n[t/\De_n]}~=~Y_0+\st\dd Y.
\eee
As soon as $Y$ is c\`adl\`ag (= right continuous with left limits), we have
$Y^{(n)}\tosc Y$ ($\om$-wise convergence for the Skorokhod
topology). If a process $Y$ belongs to the set $\vv$ of all processes
of locally of finite variation, we denote by $v(Y)_t=\int_0^t|dY_s|$
its `` variation process''.

Next we give a series of notational conventions for the convergence of
a sequence $(Y^n)$ of (c\`adl\`ag) processes; below, $\al_n$ is a
sequence of positive, possibly random, numbers\,:

\noindent $\bullet$ $Y^n\toucp Y$ or $Y^n_t\toucp Y_t$ (or,  {\em
  converges u.c.p.}) , if $\sup_{s\leq
t}|Y_s^n-Y_s|\toop 0$ for all $t>0$;

\noindent $\bullet$ $Y^n\toSp Y$ or $Y^n_t\toSp Y_t$, if the
convergence takes place in probability, for the Skorokhod topology;

\noindent $\bullet$ $Y^n\tovp Y$ or $Y^n_t\tovp Y_t$ (or, {\em converges
v.p.}) if $v(Y^n-Y)_t\toop 0$ for all $t>0$;

\noindent $\bullet$ $Y^n\tols Y$ or $Y^n_t\tols Y_t$ if there is stable
convergence in law, see below;

\noindent $\bullet$ $Y^n=\tou(\al_n)$ or $Y^n_t=\tou(\al_n)$
if $Y^n/\al_n\toucp0$;

\noindent $\bullet$ $Y^n=\tOu(\al_n)$ or $Y^n_t=\tOu(\al_n)$
if the sequences $(\sup_{s\leq t}|Y^n_s/\al_n|)_{n\geq1}$ are tight;

\noindent $\bullet$ an array $(\ze^n_i)$ of variables is {\em asymptotically
  negligible}, (AN) for short, if $\st\ze^n_i\toucp0$.

When each $Y^n$ is defined on $(\Om,\f,\PP)$, recall (see e.g. \cite{JS})
that $Y^n\tols Y$ means that $Y$ is a c\`adl\`ag process defined on an
extension of $(\Om,\f,\PP)$, and that
$\E(Zg(Y^n))\to\E(Zg(Y))$ for all bounded $\f$--measurable variable
$Z$ and all bounded continuous function $g$ on the space of all
c\`adl\`ag functions, endowed with the Skorokhod topology.

Throughout, the following functions $h_r$ for $r\in(0,\infty)$ and
$\psi_\eta$ for $\eta\in(0,\infty]$
and $\phi_s$ for $s\in[0,2]$ will often
occur : we first fix a $C^\infty$ function $\psi$ having $1_{\{|x|\leq1\}}
\leq\psi(x)\leq 1_{\{|x|\leq2\}}$, and then set
\bee\label{I3}
\left.\begin{array}{l}
h_r(x)=|x|^r,\\
\psi_\eta(x)~=~\left\{\begin{array}{ll}
\psi(x/\eta) &\mbox{if }~\eta<\infty\\
1&\mbox{if }~\eta=\infty,\end{array}\right.
\\\phi_r(x)~=~\left\{\begin{array}{ll}
1\bigwedge|x|^r&\mbox{if }~~0<r<\infty\\
1_{\R\backslash\{0\}}(x) &\mbox{if }~~r=0. \end{array}\right.
\end{array}\right\}
\eee

Next, we introduce several classes of functions on $\R$. We
denote by $\ea$ the set of all Borel functions with at most
polynomial growth, and for
$r\in(0,\infty)$ we denote by $\ea_r$ and $\ea'_r$
and $\ea''_r$ the following sets of functions :
\bee\label{I10}
\left.\begin{array}{lll}
\ea_r&:&\mbox{all $f\in\ea$ with ~$f(x)=|x|^r$~ on a neighborhood of $0$}\\
\ea'_r&:&\mbox{all $f\in\ea$ with
~$f(x)\sim|x|^r$~ as $x\to0$}\\
\ea''_r&:&\mbox{all $f$ locally bounded with ~$f(x)=$ O$(|x|^r)$~
as $x\to0$}\\
\ea'''_r&:&\mbox{all $f$ locally bounded with ~$f(x)=$ o$(|x|^r)$~
as $x\to0$}.\end{array}\right\}
\eee
We write $\ea_r^b$, $\ea'^b_r$, $\ea''^b_r$ and $\ea'''^b_r$ for the
sets of bounded functions belonging to $\ea_r$, $\ea'_r$, $\ea''_r$
and $\ea'''_r$ respectively. We have $\phi_r\in\ea^b_r\cap C^0$, where as usual
$C^p$ denotes the set of $p$ times continuously
differentiable functions, resp. continuous, for $p\geq1$, resp. $p=0$.

Below, $K$ is a constant which changes from line to line and
may depend on $X$ and its characteristics, and we write $K_p$ if
we want to emphasize its dependency on some parameter $p$. We
write $U$ for a generic $\n(0,1)$ variable, and $m_r=\E(|U|^r)$ is
its $r$th absolute moment. We also denote by $\rho_s$ the
normal law $\n(0,s^2)$, and write $\rho_s(f)=\int f(x)\rho_s(dx)$.

\subsection{The assumptions.}

We start with a semimartingale $X$ on a stochastic
basis $\proba$. We fix a truncation function $\ka$ (bounded with
compact support, with $\ka(x)=x$ on a neighborhood of $0$): this
function is a priori arbitrary and usually $\ka(x)=x1_{\{|x|\leq1\}}$,
but in this paper we suppose that it is {\em continuous} : this
simplifies some of the assumptions below. We call $(B,C,\nu)$ its
predictable characteristics : $\nu$ is the compensator of the
jump measure $\mu$ of $X$, and $C=\langle X^c,X^c\rangle$, where
$X^c$ is the continuous martingale part of $X$, and $B$ depends on the
choice of $\ka$. With $\ka'(x)=x-\ka(x)$, we then have
\bee\label{I1}
X=X_0+B+X^c+\ka\is+\ka'\star\mu.
\eee
Here and below we use
standard notation for stochastic integrals and characteristics, see for
example \cite{JS} for all unexplained notation.

We are interested in the associated processes
$V^n(f)$ and $V'^n(f)$ in (\ref{I4}) (written as $V^n(f;X)$ and
$V'^n(f;X)$ if we want to emphasize the dependency upon $X$).
For simplicity we write $\PP^n_i$ and $\E^n_i$ for the conditional
probability and expectation w.r.t. $\f_{i\De_n}$. We also
introduce some related notation, where $f$ is a small enough function
(e.g. bounded) :
\bee\label{I12}
H^n_i(f)=\E^n_{i-1}(f(\dd X)),\qquad
K^n_i(f)=\E^n_{i-1}(f(\dd X/\rn)),
\eee
\bee\label{I13}
\BH^n(f)_t:=\st H^n_i(f),\qquad \BK^n(f)_t:=\st K^n_i(f).
\eee

Our first key result needs no special assumption,
but stating it requires some additional notation: first, $C^{0,\nu}$
denotes the set of all functions on $\R$ which are
$\nu(\om;\R_+\times~dx)$--a.e. continuous, for $\PP$--almost all $\om$.
Next, we set
\bee\label{I7}
I=\{r\geq0:\phi_r\star\nu_t<\infty~~\forall t>0\}.
\eee
This is an interval of the form
$[\al,\infty)$ or $(\al,\infty)$, for some $\al\in[0,2]$. We
have $2\in I$ always, and we have $X-X^c\in\vv$ if and only if $1\in
I$, and $X$ has a.s. finitely many jumps on each finite time interval
if and only if $0\in I$. Set 
\bee\label{I8}
\left.\begin{array}{l} X'=X-X^c-X_0,\\[2mm]
1\in I\quad\Rightarrow\quad \BB=B-\ka\star\nu,\quad
X''_t=\sum_{s\leq t}\De X_s.\end{array}\right\}
\eee
So if $1\in I$ we have $X'=\BB+X''$, and $\BB$ is the ``genuine''
drift. In this case $\BB\in\vv$.

The other results need various assumptions, which we presently describe.
\vst

\nib Hypothesis (H) : \rm The characteristics $(B,C,\nu)$ of $X$
have the form
\bee\label{AS1}
B_t=\int_0^tb_sds,\qquad C_t=\int_0^tc_sds,\qquad
\nu(dt,dx)=dt~F_t(dx).
\eee
Moreover the processes $(b_t)$ and $(F_t(\phi_2))$ are locally
bounded predictable (where
$F_t(f)=\int f(x)F_t(dx)$), and the process $(c_t)$ is c\`adl\`ag adapted. \qed
\vst

Clearly (H) implies the quasi-left-continuity of $X$. Under (H), we write
\bee\label{AS2}
\si_t~=~\sqrt{c_t}.
\eee

As is well known, the form (\ref{AS1}) of the characteristics of $X$
is equivalent to the fact that $X$ can be written as
\bee\label{AS3}
X_t=X_0+\int_0^tb_sds+\int_0^t\si_sdW_s+\ka(\de)
\star(\umu-\unu)_t+\ka'(\de)\star\umu_t,
\eee
where

\noindent 1) $\si$ is given by (\ref{AS2}) and $\de$ is a ''predictable''
map from $\Om\times\R_+\times\R$ on $\R$, connected with $F_t$ by the
fact that $F_t(\om,dx)$ is the image of the Lebesgue measure on $\R$ by the map
$x\mapsto\de(\om,t,x)$.

\noindent 2) $W$ and $\umu$ are a Wiener process and a Poisson
random measure on $\R_+\times\R$ on the filtered space $\proba$
and the predictable compensator of $\umu$ is
$\unu(ds,dx)=ds\otimes dx$ (we may have to enlarge a bit the
original space in order to accommodate the pair $(W,\umu)$).
\vst

\nib Hypothesis (K) : \rm (H) holds and in (\ref{AS3}) the coefficient
$\de$ satisfies $|\de(\om,t,x)|\leq\ga_k(x)$ for all $t\leq T_k(\om)$,
where $\ga_k$ are (deterministic) functions on $\R$ with
$\int\phi_2\circ\ga_k(x)~dx<\infty$, and
$(T_k)$ is a sequence of stopping times increasing to $+\infty$.\qed
\vst

In the next hypothesis we assume that the space also supports a second
Wiener process $W'$ independent of $W$. Note that the particular form of
$\umu$ in (\ref{AS3}) or in (\ref{AS4}) below is actually irrelevant,
it could be a Poisson random measure
on $\R_+\times E$ for any space $E$ and with a compensator of the form
$dt\otimes\uF(dx)$, provided the measure $\uF$ is infinite and without
atom; or, we could have two different Poisson random measure in
(\ref{AS3}) and in (\ref{AS4}).
\vst

\nib Hypothesis (L-$s$) : \rm (H) holds and the process $\si$
in the formula (\ref{AS3}) has the form
\bee\label{AS4}
\si_t=\si_0+\int_0^t\Wb_sds+\int_0^t\Wsi_sdW_s+\int_0^t\Wsi'_sdW'_s
+\ka(\Wde)\star(\umu-\unu)_t+\ka'(\Wde)\star\umu_t,
\eee
and

a) the process $(\Wb_t)$ is optional and locally bounded;

b) the processes $(b_t)$, $(\Wsi_t)$, $(\Wsi'_t)$ are adapted
left-continuous with right limits and locally bounded;

c) the functions $\de(\om,t,x)$ and $\Wde(\om,t,x)$ are
predictable, left-continuous with right limits in $t$,
and $|\de(\om,t,x)|\leq\ga_k(x)$ and
$|\Wde(\om,t,x)|\leq\Wga_k(x)$ for all $t\leq T_k(\om)$, where $\ga_k,\wga_k$
are (deterministic) functions on $\R$ with
$\int\phi_s\circ\ga_k(x)~dx<\infty$ (with $0^0=0$) and
$\int\phi_2\circ\Wga_k(x)~dx<\infty$,
and
$(T_k)$ is a sequence of stopping times increasing to $+\infty$.\qed
\vst

In (L-$s$) we implicitly assume $s\in[0,2]$. Note that
if $s\leq s'\leq2$, then
(L-$s'$) $\Rightarrow$ (L-$s$) $\Rightarrow$ (K) $\Rightarrow$ (H),
and (L-$s$) implies that $s\in I$. It is worthwhile to emphasize that
(L-$0$) implies that $X$ has locally finitely many jumps, and also
that when $X$ is continuous then {\em all hypotheses} (L-$s$) for
$s\in[0,2]$ are identical.

Finally we have an assumption of a different nature :
\vst

\nib Hypothesis (H') : \rm We have (H) and the processes
$(c_t)$ and $(c_{t-})$ do not vanish. \qed

\begin{rem}\label{R100} \rm These assumptions, and especially (L-$s$),
may appear complicated to check. However, if $X$ is one of
the components of the solution of an SDE of the form
$d\BX_t=f(\BX_{t-})dZ_t$, where $Z$ is a multidimensional L\'evy
process and $f$ is a $C^2$ function with linear growth and locally
bounded second derivative, then (L-$2$) is automatically
satisfied. The same holds for solutions of SDEs driven by $W$ and
$\umu$. \qed
\end{rem}

\subsection{The laws of large numbers.} \label{ssLLN}

First, we have a result valid with no assumption at all on $X$ (recall the
notation (\ref{I4}), (\ref{I7}) and (\ref{I8})) :

\begin{theo}\label{T1} (i) The processes $V^n(f)$ converge in
probability in the Skorokhod sense to a suitable limit $V(f)$ in the
following cases:

\noindent(a) With $V(f)=f\star\mu$, when

\hskip2cm [a-1] $f\in\ea_2'''\cap C^{0,\nu}$,

\hskip2cm[a-2] $f\in\ea_r''\cap C^{0,\nu}$ if $r\in I\cap(1,2)$ and $C=0$,

\hskip2cm[a-3] $f\in\ea_1'''\cap C^{0,\nu}$ if $1\in I$ and $C=0$,

\hskip2cm[a-4] $f\in\ea_r''\cap C^{0,\nu}$ if $r\in I\cap(0,1]$ and
$C=\BB=0$.

\noindent(b) With $V(f)=f\star\mu+C$, when $f\in\ea'_2\cap
C^{0,\nu}$.

\noindent(c) With $V(f)=f\star\mu+v(\BB)$, when $f\in\ea'_1\cap C^{0,\nu}$
and $C=0$ and $1\in I$.

(ii) Moreover in (a) and (c) above we also have $V^n(f)-V(f)^{(n)}\tovp0$.
\end{theo}

When $f=h_r$ the case (b) ($r=2$) is well known (convergence of
the realized quadratic variation), and (a) for $r>2$ may be
found in \cite{L} for general semimartingales, and (c) ($r=1$) is also well
known because $V(f)$ is then the variation process of $X$.

The next LLNs are obtained after centering or normalization. For the
first one we need to introduce the process
\bee\label{I11}
\Si(f,\psi_\eta)=(f\psi_\eta)\is+(f(1-\psi_\eta)\star\mu,
\eee
which is well defined for $\eta\in(0,\infty]$ as soon as
$f^2\in\ea''_r$ for some $r\in I$, and also for $\eta=\infty$ if
further $f$ is bounded (it is then a locally square-integrable
martingale).

\begin{theo}\label{T2} Assume that $X$ is quasi-left-continuous.
Let $f\in\ea''_r\cap C^{0,\nu}$ for some $r\in(1,2)$. Then
  $V^n(f)-\BH^n(f\psi_\eta)\toSp\Si(f,\psi_{\eta})$ if $\eta<\infty$,
  and also if $\eta=\infty$ when $f$ is bounded.
\end{theo}

\begin{theo}\label{T3} Assume (H). Then:

(i) $\De_nV'^n(g)_t\toucp\int_0^t\rho_{\si_u}(g)du$ if $g$ is a
continuous function, in $\ea$ when $X$ is continuous, and with
$g(x)/x^2\to0$ as $|x|\to\infty$ otherwise.

(ii) $\De_n^{1-r/2}V^n(f)_t\toucp m_r\int_0^tc_u^{r/2}du$ if
$f\in\ea'_r$ and $r\in(0,2)$.

(iii) $V''^n(\vpi,\al)\toucp C_t$ for all $\vpi\in(0,\frac12)$ and $\al>0$.
\end{theo}

\begin{rem}\label{R1} \rm Theorem \ref{T2} is an LLN because the
convergence holds in probability, but it can also be viewed as a
CLT since the limiting process is a (local) martingale as soon as
$f$ is bounded and $\eta=\infty$.  \qed
\end{rem}

\begin{rem}\label{R2} \rm
Theorem \ref{T2} overlaps with (i) of
  Theorem \ref{T1}, but in the overlapping cases the two are
  of course consistent. When Theorem \ref{T2} applies and Theorem
  \ref{T1} fails, there is $t>0$ such that both sequences $(V^n(f)_t)$
  and $(\BH^n(f\vp)_t)$ are {\em not}\, tight.

When $r\in(1,2)$ and $f\in\ea'_r\cap
C^{0,\nu}$, Theorems \ref{T2} and \ref{T3}-(ii) also overlap:
an equivalent way of writing the later is
$\De_n^{1-r/2}\left(V^n(f)-\BH^n(f)\right)\toucp 0$ (see the proofs
below), so Theorem \ref{T2} in this
case is the CLT associated with the LLN of Theorem \ref{T3}-(ii)
in a sense, or perhaps rather as a ''second order'' LLN because the
convergence takes place in probability. \qed
\end{rem}

\begin{rem}\label{R3} \rm The reader will note the - different -
assumptions in the last two theorems. Theorem \ref{T2} probably fails
if $X$ is not quasi-left continuous. Theorem \ref{T3} just makes no
sense if (H) fails (or rather, if the second equality in
(\ref{AS1}) fails), and the quasi-left continuity is by no
means enough for it. \qed
\end{rem}

\subsection{The central limit theorems.}

The various CLTs below involve stable convergence in law,
for which we need some ingredients.
Consider an auxiliary space $(\Om',\f',\PP')$ supporting
a $d$-dimensional Brownian motion $\BW=(\BW^j)_{1\leq j\leq d}$, two
sequences $(U_n)$ and $(U'_n)$ of
$\n(0,1)$ variables, and a sequence $(\ka_n)$ of variables uniformly
distributed on $(0,1)$, all of these being mutually independent. Then put
\bee\label{C1}
\WOm=\Om\times\Om',\qquad \Wf=\f\otimes\f',\qquad\WP=\PP\otimes\PP'.
\eee
and extend the variables $X_t,~b_t$, ... defined on $\Om$ and
$\BW,~U_n$,... defined on $\Om'$ to the product $\WOm$ in the obvious way,
without changing the notation. We write $\WE$ for the expectation
w.r.t. $\WP$. Finally, denote by $(T_n)_{n\geq1}$ an enumeration of
the jump times of $X$ which are stopping times, and let $(\Wf_t)$ be
the smallest (right-continuous) filtration of $\Wf$ containing the
filtration $(\f_t)$ and w.r.t.\ which $\BW$ is
adapted and such that $U_n$ and $U'_n$ and $\ka_n$ are
$\Wf_{T_n}$-measurable for all $n$.

Obviously, $\BW$ is an $(\Wf_t)$-Brownian motion under $\WP$,
as well as $W$, and $W'$ under (L-$2$), whereas $\umu$ is still a Poisson
measure with compensator $\unu$ for this bigger filtration. The
dimension $d$ of $\BW$ is the number of processes for which we want to
have a joint CLT in Theorem \ref{T8} below, in the other theorems we
have $d=1$ and we then write $\BW^1=\BW$.

The limiting processes we obtain below are of the form $Y=(Y^j)_{1\leq
j\leq d}$ with $Y^j_t=\sum_{k=1}^d\int_0^t\te^{jk}_u~d\BW^k_u$ for
suitable $(\f_t)$-adapted $d\times d$-dimensional
c\`adl\`ag processes $(\te_t)$, or the sum of $Y_t$ plus a
process of the form
\bee\label{C2}
Z(g)_t=\sum_{p:~T_p\leq t} g(\De X_{T_p})\Big(\sqrt{\ka_p}~U_p~\si_{T_p-}
+\sqrt{1-\ka_p}~U'_p~\si_{T_p}\Big),
\eee
for some function $g\in\ea''_1$. As we will
check in Lemma \ref{LC4} below, this formula defines a semimartingale on
the extended space, whose conditional law w.r.t.\ $\f$ depends on the
processes $X$ and $c$ (or $\si$) but not on the particular choice of
the stopping times $T_n$. Moreover, again conditionally on $\f$, the
two processes $Y$ and $Z(g)$ are independent and are {\em martingales}
  with variance-covariance given by
\bee\label{C2'}
\left.\begin{array}{l}
\WE(Y^j_tY^k_t\mid\f)~=~\int(\te_u\te^\star_u)^{jk}~du\\[2mm]
\WE(Z(g)^2_t\mid\f)~=~C(g)_t~:=~\sum_{p:~T_p\leq t}g(\De X_{T_p})^2
(c_{T_p-}+\frac12~\De c_{T_p}),\end{array}\right\}
\eee
where $\te^\star$ is the transpose and $\De c_{T_p}$ is the jump of
the process $(c_t)$ at time
$T_p$. Moreover, conditionally on $\f$,  $Y$ is even a {\it Gaussian
  martingale},  and $Z(g)$ also as soon as the processes
$X$ and $\si$ have no common jumps. This will also be checked later.

Now we state a variety of CLTs, related with some of the LLNs given above,
although the picture is far from being complete. As said before, Theorem
\ref{T2} is already a CLT in a sense, and we start with a result
extending this theorem to the case $r=1$. The
other CLTs are related to Theorem \ref{T3} and with a special case of
Theorem \ref{T1}-(a), with unfortunately some
unwanted restrictions. We complement these CLTs with some ``tightness''
results, in view of applications. Finally we will end up with a
multidimensional CLT which contains the previous results and is
complicated to state, but which probably is the most useful result for
practical applications, at least for those we have in mind.

\begin{theo}\label{T4} Assume (H). Let $f\in\ea'_1\cap C^{0,\nu}$ and
$\eta<\infty$, or $\eta=\infty$ if $f$ is bounded. Then
$V^n(f)-\BH^n(f\psi_\eta)\tols\Si(f,\psi_\eta)_t+
\sqrt{m_2-m_1^2}\int_0^t\si_u~d\BW_u$~~ \rm (note that $m_2-m_1^2=1-2/\pi$).
\end{theo}

\begin{theo}\label{T5} Assume (L-$s$), and let $g$ be an even $C^2_b$ 
  function on $\R$.

 (i) ~$\frac1{\rn}\left(\De_nV'^n(g)_t-\int_0^t\rho_{\si_u}(g)du\right)\!\!
\tols \!\!\int_0^t\sqrt{\rho_{\si_u}(g^2)-(\rho_{\si_u}(g))^2}~d\BW_u$
if $s\leq1$;

(ii)~ $\De_nV'^n(g)_t-\int_0^t\rho_{\si_u}(g)du=\tou(\De_n^{1-s/2})$
otherwise. 

\ni When $X$ is continuous, we have (i) under (L-$2$), as soon as $g$
is $C^1$ and even and $g'\in\ea$. 
\end{theo}

\begin{theo}\label{T6} Assume (L-$s$) and (H'), and let
$f\in\ea_r$ for some $r\in(0,1]$.

 (i) ~$\frac1{\rn}\left(\De_n^{1-r/2}V^n(f)_t-m_r\int_0^tc_u^{r/2}
\right)\!\!\tols\!\! \sqrt{m_{2r}-m_r^2}\int_0^tc_u^{r/2}~d\BW_u$ if either
$s\leq\frac23$ and $r<1$, 
or $\frac23<s<1$ and $\frac{1-\sqrt{3s^2-8s+5}}{2-s}<r<1$;

 (ii) ~$\De_n^{1-r/2}V^n(f)_t-m_r\int_0^tc_u^{r/2}=
\tou\Big(\De_n^{\frac{(2-s)(1+r)(2-r)}{4+2s(1-r)}-\ep}\Big)$ for all
$\ep>0$, otherwise. 

\ni When $X$ is continuous, we have (i) under (L-$2$) and (H') when
$r\in(0,1]$, and also under (L-$2$) only when $r>1$.
\end{theo}

\begin{theo}\label{T6'} Assume (L-$s$), and let $\vpi\in(0,\frac12)$
  and $\al>0$. Then

 (i) ~$\frac1{\rn}\left(V''^n(\vpi,\al)_t-C_t
\right)\!\!\tols \!\!\sqrt{2}\int_0^tc_u~d\BW_u$ if
$s\leq\frac{4\vpi-1}{2\vpi}$ (hence $\vpi\geq\frac14$ and $s<1$); 

 (ii) ~$V''^n(\vpi,\al)_t-C_t=\tou(\De_n^{(2-s)\vpi})$ otherwise.
\end{theo}

\begin{theo}\label{T7} Let $f$ be a $C^1$ function on $\R$.

(i) Under (K), and if $f$ is $C^2$ on a neighborhood of $0$, with
$f(0)=f'(0)=0$ and $f''(x)=$ o$(|x|)$ as $x\to0$, then
$\frac1{\sqrt{\De_n}}~(V^n(f)_t-V(f)^{(n)}_t)\tols Z(f')_t$
(with $V(f)=f\star\mu$).

(ii) Under (L-$2$) and if $f\in\ea_2$ we
have $\frac1{\sqrt{\De_n}}~(V^n(f)_t-V(f)^{(n)}_t)\tols Z(f')_t
+\sqrt{2}\int_0^tc_u~d\BW_u$ (with $V(f)=C+f\star\mu$).
\end{theo}

\begin{rem}\label{R4} \rm We do not have stable convergence in law of
the processes $\frac1{\sqrt{\De_n}}~(V^n(f)-V(f)_t)$ in the last
theorem, and not even mere convergence in law, because of some
peculiarity of the Skorokhod topology. However these processes
converge finite-dimensionally stably in law to the limits described
above. \qed
\end{rem}

\begin{rem}\label{R6} \rm The limiting process in (ii) of Theorem
\ref{T7} looks pretty much like the limiting process obtained in
\cite{JP} for the error term in the Euler approximation of the
solution of SDEs driven by L\'evy processes. This is of course not
just by chance ! \qed 
\end{rem}

\begin{rem}\label{R5} \rm Suppose that $f=h_r$. We have a CLT
for $V^n(f)$ in the following cases~:

$\bullet$ if $r<1$, in Theorem \ref{T5} (subject to some - perhaps
unnecessary - restrictions on the value of $s$ for which (L-$s$) holds), after
normalization and centering;

$\bullet$ if $r=1$, in Theorem \ref{T4}, after centering;

$\bullet$ if $1<r<2$, in Theorem \ref{T2}, after centering;

$\bullet$ if $r=2$ or $r>3$, in Theorem \ref{T7}, after
normalization and centering.

\noindent When $2<r\leq3$, there is no CLT, at least with the
natural centering of the associated LLN, although a CLT with a more
adequate centering might exist: see \cite{J2} for a more thorough
description of this fact when $X$ is a L\'evy process.  \qed
\end{rem}

Finally, we give the announced multidimensional CLT, in
which we consider components as in Theorems \ref{T5}, \ref{T6}, \ref{T6'} and
\ref{T7}. Below we have a $d$--dimensional process and the index set
$\{1,\ldots,d\}$ for the components is partitioned into five
(possibly empty) subsets $J_l$. We consider the
process $Y^n=(Y^{n,j})_{1\leq j\leq d}$ having the following components :
\bean
\bullet j\in J_1~&\Rightarrow~& Y^{n,j}_t=\De_nV'^n(f_j)_t-
\int_0^t\rho_{\si_u}(f_j)du,~\mbox{where $f_j$ is $C^2_b$ and even;} \\
\bullet j\in J_2~&\Rightarrow~& Y^{n,j}_t=\De_n^{1-r(j)/2}V^n(f_j)_t-
m_{r(j)}\int_0^t\si_u^{r(j)/2}du,~\mbox{ where $f_j\in\ea_{r(j)}$}\\
&&~~\mbox{for some $r(j)\in(0,1)$ in general or $r(j)\in(0,\infty)$ if
  $X$ is continuous;}\\ 
\bullet j\in J_3~&\Rightarrow~&Y^{n,j}_t=V''^n(\vpi_j,\al_j)_t-C_t,\\
&&~~\mbox{where $\vp_j\in[1/4,1/2)$ and $\al_j>0$; we then put $r(j)=2$;}\\
\bullet j\in J_4~&\Rightarrow~&Y^{n,j}=V^n(f_j)-V(f_j)^{(n)},~~
\mbox{ where
$f_j$ is $C^1$, and $C^2$ on a neighbor-}\\
&&~~\mbox{hood of $0$ with $f_j(0)=f'_j(0)=0$ and
$f''_j(x)=$o$(|x|)$ as $x\to0$;}\\
\bullet j\in J_5~&\Rightarrow~&Y^{n,j}=V^n(f_j)-V(f_j)^{(n)},~~
\mbox{ where $f_j\in\ea_2\cap C^1$; we then put $r(j)=2$.}
\eean

\begin{theo}\label{T8} With the previous setting, we assume (H') if
$J_2\neq\emptyset$, and (L-$s$) for some $s\in[0,2]$ satisfying 

$J_1\neq\emptyset~~\Rightarrow~~s<1$, 

$J_2\neq\emptyset~~\Rightarrow~~$ either $0\leq s\leq\frac23$ or
$\frac23<s<1$ and 
$\frac{1-\sqrt{3s^2-8s+5}}{2-s}<\inf_{j\in J_2}r(j)$, 

$J_3\neq\emptyset~~\Rightarrow~~s<\inf_{j\in J_3}\frac{4\vpi_j-1}{2\vpi_j}$.

\noindent Then $\frac1{\rn}~Y^n~\tols Y$, where
\bee\label{C100}
Y^j_t=\left\{\begin{array}{ll}
\sum_{k\in J_1\cup J_2\cup J_3\cup J_5}\int_0^t\te^{jk}_u~d\BW^k_u\qquad&
\mbox{if }~j\in J_1\cup J_2\cup J_3\\[2mm]
Z(f'_j)_t\qquad&\mbox{if }~j\in J_4,\\[2mm]
Z(f'_j)_t+\sum_{k\in J_1\cup J_2\cup J_3\cup J_5}\int_0^t\te^{jk}_u~
d\BW^k_u\qquad&\mbox{if }~j\in J_5,\end{array}\right.
\eee
and where $\te=(\te^{jk})_{j,k\in J_1\cup J_2\cup J_4}$ is an
$(\f_t)$-adapted c\`adl\`ag process whose square $\te\te^\star$
is the symmetric matrix characterized by
\bee\label{C200}
(\te_t\te_t^\star)^{jk}=\left\{\begin{array}{ll}
\rho_{\si_t}(f_jf_k)-\rho_{\si_t}(f_j)\rho_{\si_t}(f_k),
\quad&j,k\in J_1\\[2mm]
(m_{r(j)+r(k)}-m_{r(j)}m_{r(k)})c_t^{r(j)/2+r(k)/2},
&j,k\in J_2\cup J_3\cup J_5\\[2mm]
\rho_{\si_t}(h_{r(j)}f_k)-\rho_{\si_t}(h_{r(j)})\rho_{\si_t}(f_k),
\quad& j\in J_2\cup J_3\cup J_5,~k\in J_1.
\end{array}\right.
\eee

\ni When $X$ is continuous, the same holds under (L-$2$) and (H') as
soon as the $f_j$'s for $j\in J_1$ are $C^1$ and even with $f_j'\in\ea$,
and $r(j)\in(0,\infty)$ for $j\in J_2$, and one can relax (H') if
$r(j)>1$ for all $j\in J_2$.
\end{theo}

(It is easy to check that the right side of (\ref{C200}) is a positive
symmetric matrix indexed by $J_1\cup J_2\cup J_3\cup J_5$, so it has a
``square-root'' $\te_t$).

\section{Theorems \ref{T1} and \ref{T2}}\label{secKR}
\setcounter{equation}{0}
\renewcommand{\theequation}{\thesection.\arabic{equation}}

\subsection{Proof of Theorem \ref{T1}.}

The idea of the proof is the same as in \cite{J2}, but the details are
slightly more involved, so we give a complete proof.

\ni\it Step 1. \rm If $f$
satisfies any one of the conditions in (i) the process $f\star\nu$
is in $\vv$, hence $V(f)=f\star\mu$ as well. In view of the convergence
$V(f)^{(n)}\tosc V(f)$ it is clear that (ii)
implies (i-a) and (i-c). Below, we use the notation
\bee\label{KR0}
Z^n(f)=V^n(f)-V(f)^{(n)}.
\eee

\ni\it Step 2: \rm Here we prove (i) and (ii) when
$f\in C^{0,\nu}$ vanishes on a neighborhood of $0$, say
$[-2\ep,2\ep]$, hence $V(f)=f\star\mu$. For any fixed $\ep>0$ we set :
\bee\label{KR1}
\left.\begin{array}{l}
\bullet~ S_1,S_2,\cdots~~\mbox{are the successive jump times of $X$ with
}~~|\De X_t|>\ep,\\
\bullet~ R_p=\De X_{S_p},\\
\bullet~ X(\ep)_t=X_t-(x1_{\{|x|>\ep\}})*\mu_t=X_t-\sum_{p:~S_p\leq t}R_p,\\
\bullet~ R'^n_p=\dd X(\ep)~~ \mbox{on the set }~
\{(i-1)\De_n<S_p\leq i\De_n\},\\
\bullet~ \Om_n(T,\ep)~~\mbox{is the set of all $\om$ such that each 
interval}~[0,T]\cap((i-1)\De_n,i\De_n]\\
\quad~\mbox{contains at most one $S_p(\om)$, and
that $|\dd X(\ep)(\om)|\leq2\ep$ for all $i\leq T/\De_n$.}
\end{array}\right\}
\eee
All these depend on $\ep$ of course, and $\Om_n(T,\ep)\to\Om$ as
$n\to\infty$.

Recalling $f(x)=0$ when $|x|\leq2\ep$, we see that on the
set $\Om_n(T,\ep)$ and for all $t\leq T$,
\bee\label{KR2}
v(Z^n(f))_t=\sum_{p:~S_p\leq\De_n[t/\De_n]}
|(f(R_p+R'^n_p)-f(R_p))|,
\eee
Since $f\in C^{0,\nu}$ there
is a null set $N$ such that, if $\om\notin N$, then $f$ is
continuous at each point $R_p(\om)$, whereas $R'^n_p(\om)\to0$, so
$v(Z^n(f))_T\to0$ when $\om\notin N$. Hence
(ii) is obvious (we even have almost sure convergence).
\vst

\ni\it Step 3: \rm Here we prove (ii) in case (c), so we assume $1\in
D$ and $C=0$. As said before, $X\in\vv$ and $v(X-X_0)=V(h_1)$ (recall
(\ref{I3})), and it
is well known that $V^n(h_1)_t$ converges pointwise to $V(h_1)_t$. Then
$Z^n(h_1)_t\to 0$ and, since $Z^n(h_1)$ is a nonpositive decreasing
process, we in fact have $v(Z^n(h_1))_t\to0$ for all $t$.

Now let $f\in\ea'_1\cap C^{0,\nu}$. We have $|(f-h_1)\psi_\eta|\leq
\ep_\eta h_1$, where $\ep_\eta\to0$ as $\eta\to0$. We have
$$v(Z^n(f))~\leq~(1+\ep_\eta)v(Z^n(h_1))+v(Z^n((f-h_1)(1-\psi_\eta)))
+\ep_\eta V(h_1).$$
The first two terms on the right go to $0$ a.s.\ by the above and Step
2, and $V(h_1)$ is finite-valued and $\ep_\eta\to0$, hence the result.
\vst

\ni\it Step 4: \rm Here we prove the remaining claims (ii-a)
and (i--b), assuming that
\bee\label{KR3}
Z^n(h_r\psi_\eta)\toucp0
\eee
in the relevant cases: that is either $r=2$ (hence $V(h_2\psi_\eta)=
C+(h_2\psi_\eta)\star\mu$), or $r\in I\cap(1,2)$ and $C=0$, or $r\in
I\cap(0,1)$ and $C=\BB=0$ (so $V(h_2\psi_\eta)=(h_2\psi_\eta)\star\mu$ in
these two cases).

Assume $f\in\ea_2'''\cap C^{0,\nu}$. Then $|f\psi_\eta|\leq
\ep_\eta h_2\psi_\eta$, with $\ep_\eta\to0$ as $\eta\to0$, and thus
$$v(Z^n(f))\leq v(Z^n(f(1-\psi_\eta))
+\ep_\eta(|Z^n(h_2\psi_\eta)|+2V(h_2\psi_\eta)).$$
The first term on the right goes to $0$ a.s. by Step 2, so (\ref{KR3})
and $\ep_\eta\to0$ and $V(h_2\psi_\eta)_t<\infty$ give (ii) in case
[a-1]. When $f\in\ea_1'''\cap
C^{0,\nu}$ and $1\in I$ and $C=0$, the same argument with $h_1$
instead of $h_2$ works (use Step 3), and we have (ii) in case [a-3].

When $f\in\ea_r''\cap C^{0,\nu}$ with $r<2$, we have $|f\psi_\eta|\leq
Kh_r\psi_\eta$ for all $\eta$ small enough, hence
$$v(Z^n(f))\leq
v(Z^n(f(1-\psi_\eta)))+K|Z^n(f_r\psi_\eta)|+2K(f_r\psi_\eta)\star\mu  .$$
The first two terms on the right go to $0$ in probability by Step 2
and (\ref{KR3}), and the third term goes to $0$ as $\eta\to0$
because $r\in I$. So we have (ii) in cases [a-2] and [a-4].

Finally let $f\in\ea'_2\cap C^{0,\nu}$, so $|(f-h_2)\psi_\eta|\leq
\ep_\eta h_2\psi_\eta$, with $\ep_\eta\to0$ as $\eta\to0$, and thus
$$|Z^n(f)|~\leq~
v(Z^n(f(1-\psi_\eta)))+(1+\ep_\eta)|Z^n(h_2\psi_\eta)|
+\ep_\eta +V(h_2),$$
and we conclude (i--b) as above.
\vst

\ni\it Step 5: \rm We are left to prove (\ref{KR3}). In other words, it
is enough to prove that if $f$ is $C^2$ outside $0$, with compact
support and $f(x)=|x|^r$ around $0$, and when either
$r=2$, or $1<r<2$ and $C=0$, or $0<r<1$ and $C=\BB=0$, then we have
$Z^n(f)\toucp0$. Set
$$g(x,y)=f(x+y)-f(x)-f(y)-\ka(x)f'(y),\quad
k(x,y)=f(x+y)-f(x)-f(y)$$
with the convention $f'(0)=0$ if $r<1$ (otherwise, $f'(0)$ is the
derivative of $f$ at $0$, of course). Recall that $V(f)=f\star\mu$
when $r<2$ and $V(f)=C+f\star\mu$ if $r=2$.

Define the process $Y^n$ by $Y^n_t=X_t-X_{(i-1)\De_n}$ for
$t\in[(i-1)\De_n,i\De_n]$. It\^o's formula when $r=2$ and
its extension as given in Theorem 3.1 of \cite{JJM} when $r<2$ give us
$$ Z^n(f)_t=
\sum_{i=1}^{[t/\De_n]}\left(f(Y^n_{i\De_n})-\dd V(f)\right)
=\sum_{i=1}^{[t/\De_n]}(A^n_i+M^n_i),$$
where (recall that $C=X^c=0$ when $r<2$ here, so $f''$ does not occur
below in that case)
\bee\label{KR4}
\left.\begin{array}{l}
A^n_i=\left\{\begin{array}{ll}
\itai\left(f'(Y^n_{s-})dB_s+(\frac12~f''(Y^n_s)-1)dC_s\right)&\\
\hskip3cm+\itai\int_\R g(x,Y^n_{s-})\nu(ds,dx)
\quad&\mbox{if }~1<r\leq2\\[2mm]
\itai\int_\R k(x,Y^n_{s-})\nu(ds,dx)
&\mbox{if }~0<r<1,
\end{array}\right.\\[1.2cm]
M^n_i=
\itai f'(Y^n_{s-})dX^c_s+
\itai\int_\R k(x,Y^n_{s-})(\mu-\nu)(ds,dx).
\end{array}\right\}
\eee
In other words, $Z^n(f)=A(n)^{(n)}+M(n)^{(n)}$, where
\bee\label{KR5}
A(n)_t=\left\{\begin{array}{ll}
\int_0^t\left(f'(Y^n_{s-})dB_s+(\frac12~ f''(Y^n_-)-1)dC_s\right)
+g(x,Y^n_-)\star\nu_t
&\mbox{if }~1<r\leq2\\[2mm]
k(x,Y^n_-)\star\nu_t&\mbox{if }~0<r<1,
\end{array}\right.
\eee
$$M(n)_t= \int_0^tf'(Y^n_{s-})dX^c_s+k(x,Y^n_-)\is_t.$$
In particular $M(n)$ is a locally square--integrable martingale, whose
predictable bracket $\langle M(n),M(n)\rangle$ is such that
$A'(n)-\langle M(n),M(n)\rangle$ is non--decreasing
(see Theorem II.1.33 of \cite{JS}), where
\bee\label{KR7}
A'(n)=f'(Y^n)^2\bullet C+k(x,Y^n_-)^2\star\nu,
\eee
\vst

\ni\it Step 6: \rm At this stage, it remains to prove that
$A(n)\toucp0$ and $M(n)\toucp0$, and for the last property Lenglart
domination property (Lemma I.3.30 of \cite{JS}) it is enough to prove
$A'(n)\toucp0$.

Suppose first that $1<r\leq2$, so the function $f$ is $C^1_b$ and $f'$ is
H\"older with index $r-1$, and $f(0)=f'(0)=0$, hence $|k(x,y)|\leq
C\phi_1(x)$ and $|g(x,y)|\leq C\phi_r(x)$, and obviously
$f'(y)$ and $k(x,y)$ and $g(x,y)$ all go to $0$ as $y\to0$ Moreover if
$r=2$ we also have $\frac12 f''(y)-1\to0$ as well. By the
assumption that $r\in I$ we have
$\phi_r\star\nu_t<\infty$, and a fortiori $\phi_1^2\star\nu_t<\infty$,
for all $t>0$. Since $Y^n_{s-}\to0$ as $n\to\infty$, we deduce from
the dominated convergence theorem and also from the property $C=0$
when $r<2$ that $\sup_{s\leq t}|A(n)_s|$ and
$\sup_{s\leq t}A'(n)_s$ both go to $0$ pointwise, and the result is proved.

Second, assume that $r<1$. Then $|k(x,y)|\leq
C\phi_r(x)$, and again $k(x,y)\to0$ as $y\to0$. Then we conclude as
above. \qed

\subsection{Some consequences.}

Now we derive some ``technical'' consequences of this basic result.

\begin{lem}\label{Lem1} Suppose that the pair $(X,f)$ satisfies one of
  the conditions of Theorem \ref{T1}--(i), and also that
  $f$ is bounded. Then we have:

(i) $\BH^n(f)=\tOu(1)$.

(ii) If $X$ is quasi--left continuous, $\BH^n(f)\toucp\BH(f)$, where
$\BH(f)$ is the (continuous) predictable compensator of $V(f)$, that is
$\BH(f)=f\star\nu$ in case (a), and $\BH(f)=C+f\star\nu$ in case (b),
and $\BH(f)=v(\BB)+f\star\nu$ in case (c).
\end{lem}

Our conditions imply that $|f|\star\mu$ is finite with
bounded jumps, so$f\star\nu$ is well defined.
We cannot hope for (ii) to be true if $X$ is not quasi--left
continuous. In general, $\BH^n(f)_t$ goes to $\BH(f)_t$
for any $t$ which is not a fixed time of discontinuity of $X$, but the
convergence is for the weak $\si(\LL^1,\LL^\infty)$ topology on
$\LL^1$: so it is not likely to be really useful~!
\vst

\nib Proof. \rm First we observe that if $f$ satisfies the assumptions
of case (a) of Theorem \ref{T1}, then the functions $f^+$, $f^-$ and
$|f|$ satisfy the same; when $f$ satisfies the assumptions of cases
(b) or (c), then $f^+$ and $|f|$ satisfy the same, whereas $f^-\in
\ea'''_2\cap C^{0,\nu}$. So it is enough to prove the result
when $f\geq 0$.

Under our assumptions, the increasing processes $V(f)$, $\BH(f)$,
$v(B)$, $C$ and $\phi_2\star\nu$ are locally bounded. So there is a
sequence $T_p$ of stopping times increasing to infinity, such that we
have identically 
\bee\label{KR10}
V(f)_{T_p}+\BH(f)_{T_p}+v(B)_{T_p}+C_{T_p}+\phi_2\star\nu_{T_p}+\leq
K_p.
\eee

Set $\BH^{n,p}(f)_t=\sum_{i=1}^{[t/\De_n]}H^{n,p}_i(h)$ and
$H^{n,p}_i(f)=\E^n_{i-1}(f(\dd X^{T_p}))$. We have
\bean
&&\E\left(\sup_{s\leq t}|\BH^n(f)_s-\BH^{n,p}(f)_s|~1_{\{T_p>t\}}\right)\\
&&\qquad\leq
\E\left(\st 1_{\{T_p>(i-1)\De_n\}}
\E^n_{i-1}\left(|f(\dd X)-f(\dd X^{T_p})|\right)\right)\\
&&\qquad\leq
K\E\left(\st 1_{\{T_p>(i-1)\De_n\}}~
\PP^n_{i-1}(T_p\leq i\De_n)\right)~\leq~K\PP(T_p\leq t),
\eean
where the second inequality above follows from $0\leq f\leq K$. Hence we
readily deduce the following
implications from the fact that $\PP(T_p\leq t)\to0$ as $p\to\infty$
for all $t$:
\bee\label{KR11}
\left.\begin{array}{ll}
\BH^{n,p}(f)_t=\tOu(1),\quad\forall p\qquad
&\Rightarrow\quad\BH^{n}(f)_t=\tOu(1),\\[2mm]
\BH^{n,p}(f)_t~\toucp~ \BH(f)_{t\bigwedge T_p},\quad\forall p\qquad
&\Rightarrow\quad\BH^{n}(f)_t~\toucp~ \BH(f)_t.
\end{array}\right\}
\eee
Therefore for (i) (resp. (ii)) it is enough to prove the first
(resp. second) left side property in (\ref{KR11}). Equivalently, it is
enough to prove the results when $X$ is such that (\ref{KR10}) holds for
$T_1=\infty$. So we proceed to proving (i) and (ii) under
this additional assumption.

(i) Set $S_{n,q}=\inf(t:V^n(f)_t\geq q)$, hence
$$\E(\BH^n(f)_{S_{n,q}})=\E(V^n(f)_{S_{n,q}})\leq q+K_1.$$
Now, $V^n(f)_t\toop V(f)_t<\infty$ for all $t$, hence
\bee\label{KR12}
\lim_{q\to\infty}~\sup_n~\PP(S_{n,q}<t)=0.
\eee
Combining these two properties gives the tightness of each
sequence $(\BH^n(f)_t)_n$.

(ii) Recall the following property, known as
the ``approximated Laplacians'' property, holds because $V(f)$ is
quasi--left continuous, see e.g. \cite{M}):
\bee\label{KR13}
\BH'^n(f)_t:=\st\E^n_{i-1}(\dd V(f))~\toucp~\BH(f)_t
\eee
(in \cite{M} the convergence is for each $t$, but it is also
u.c.p. because both sides are increasing in $t$, and $\BH(f)$ is
continuous).

Now we prove the result in cases (a) and (c) of Theorem
\ref{T1}. By (ii) of this theorem we know that
$v(V^n(f)-V(f)^{(n)})_t\toop0$ for all
$t$. By hypothesis $V(f)_\infty\leq K$, so if $S_{n,q}$ is like in
(i), we have $\E(v(V^n(f)-V(f)^{(n)})_{t\bigwedge S_{n,q}})\to0$, and
a fortiori $\E(v(\BH^n(f)-\BH'^n(f))_{t\bigwedge
S_{n,q}})\to0$. Since (\ref{KR12}) holds, we deduce the result from
(\ref{KR13}).

Finally we prove the result in case (b). Using the notation of Step 5 of
the proof of Theorem \ref{T1}, we have
$$H^n_i(f)=\dd \BH^n(f)+\E^n_{i-1}(A^n_i).$$
Then in view of (\ref{KR13}) it is enough to have $\E(v(A(n)_\infty)\to0$
(recall (\ref{KR5}, here $r=2$). But since $v(B)_\infty$, $C_\infty$,
and $\phi_2\star\nu_\infty$ are bounded, this is proved exactly as
$A(n)\toucp0$ in Step 6 of the proof of Theorem \ref{T1}. \qed
\vsc

\begin{lem}\label{Lem2} Assume $C=0$, and let $s\in I\cap[0,2]$. Let
  $f\in\ea''^b_r$ for some $r>0$. Then
\bee\label{KR14}
\BH^n(f)=\left\{\begin{array}{ll}
\tou(\De_n^{r/s-1})\qquad&\mbox{if }~r<s,~s>1\\[1.5mm]
\tOu(\De_n^{r-1})\qquad&\mbox{if }~r<1,~s\leq1\\[1.5mm]
\tOu(1)\qquad&\mbox{if }~r\geq s\bigvee1.\end{array}\right.
\eee
\end{lem}

\nib Proof. \rm There is a function $f_r\in\ea_r^b\cap C^0$ such
that $|f|\leq f_r$. Since
$|\BH^n(f)|\leq\BH^n(f_r)$ it is enough to prove the
result for $f_r$. Set $s'=s\bigvee1$, which is in $I\cap[1,2]$.

When $r\geq s'$ we have $f_r\in\ea'''_2\cap C^0$ if $r>2$, and
$f_r\in\ea_2\cap C^0$ if $r=2$, and $f_r\in\ea''_{s'}\cap C^0$ if
$1<r<2$, and $f_r\in\ea'_1\cap C^0$ if $r=1$, and $r\in I$ always, so
$f_r$ is always in one of the cases of Theorem \ref{T1}, and
the result follows from Lemma \ref{Lem1}.

When $r<s'$, H\"older inequality yields for all $\ep>0$:
$$\De_n^{1-r/s'}\BH^n(f_r)_t\leq
t^{1-r/s'}\left(\BH^n((f_r\psi_\ep))^{s'/r})_t
\right)^{r/s'}+\De_n^{1-r/s'}\BH^n(f_r(1-\psi_\ep))_t.$$
Since $f_r(1-\psi_\ep)\in\ea'''^b_2\cap C^0$, by Lemma
\ref{Lem1} again the last term above goes to
$0$ in probability for any $\ep>0$ because $r<s$. Since
$(f_r\psi_\ep)^{s'/r}\in\ea_{s'}'^b\cap C^0$ we deduce as above, from
Lemma \ref{Lem1} again, that the first
term on the right goes to $t^{1-r/s'}(f_r\psi_\ep)^{s'/r}\star\nu_t$
if $s'>1$ and to $t^{1-r/s'}\left((f_r\psi_\ep)^{s'/r}\star\nu_t
+v(\BB)_t\right)$ if $s'=1$. Now, $(f_r\psi_\ep)^{s'/r}\star\nu_t
\to0$ as $\ep\to0$ because $s'\in I$. Then we obtain the first and
second properties in (\ref{KR14}). \qed

\subsection{Proof of Theorem \ref{T2}.}

Let $f\in\ea''_r\cap C^{0,\nu}$
with $r\in(1,2)$. Let $\eta\in(0,\infty)$, or $\eta=\infty$ when
$f$ is bounded: in all cases the process $\Si(f,\psi_\eta)$ is well defined.

For any $\ep>0$ we have $f-f\psi_\ep\in\ea'''_2\cap C^{0,\nu}$, so
$V^n(f-f\psi_\ep)\toSp (f-f\psi_\ep)\star\mu$ by Theorem \ref{T1} and
$\BH^n(f\psi_\eta-f\psi_\ep)\toucp (f(\psi_\eta-\psi_\ep))\star\nu$ by
Lemma \ref{Lem1}. Therefore, as soon as $\ep<\eta$,
$V^n(f(1-\psi_\ep)-\BH^n(f(1-\psi_\ep)\psi_\eta)\toSp
\Si(f(1-\psi_\ep)),\psi_\eta)$. Moreover it is obvious that
$\Si(f(1-\psi_\ep),\psi_\eta)\toucp\Si(f,\psi_\eta)$ as
$\ep\to0$. Hence in order to prove the result it is enough to show
that if $M^n(\ep)=V^n(f\psi_\ep)-\BH^n(f\psi_\ep)$, then we have
\bee\label{N4}
t>0,~\rho>0\quad\Rightarrow\quad
\lim_{\ep\to0}~\limsup_n~\PP(\sup_{s\leq t}|M^n(\ep)|>\rho)=0.
\eee

Now the process $M^n(\ep)$ is a locally bounded martingale w.r.t. the
filtration $(\f^n_t=\f_{\De_n[t/\De_n]})_{t\geq0}$, and its
predictable quadratic variation is
$$C^n(\ep)_t=\st\left(H^n_i(f\psi_\ep)^2)-(H^n_i(f\psi_\ep))^2\right)
\leq \BH^n((f\psi_\ep)^2)_t.$$
Observe that $(f\psi_\ep)^2\in\ea_{2r}'''\cap
C^{0,\nu}$, whereas $2r>2$. Then $\BH^n((f\psi_\ep)^2)\toucp
(f\psi_\ep)^2\star\nu$ by Lemma \ref{Lem1}, hence
\bee\label{N5}
t>0,~\rho>0\quad\Rightarrow\quad
\lim_{\ep\to0}~\limsup_n~\PP(C^n(\ep)_t>\rho)=0.
\eee
By Lenglart inequality it is well known that (\ref{N5}) implies
(\ref{N4}), and we are done.

\section{Theorem \ref{T3}}\label{secCLTA}
\setcounter{equation}{0}
\renewcommand{\theequation}{\thesection.\arabic{equation}}

\subsection{Technical consequences of (H).}

The assumption (H) is ''local'', in the sense that it does not
require any integrability assumptions (in $\om$) on the characteristics.
However having ``locally bounded'' replaced by ``bounded'', for
example, simplifies a lot of technical problems. This is why we
introduce ``global'' and apparently much stronger conditions:
\vst

\nib Hypothesis (SH): \rm We have (H), and the processes
$(b_t)$, $(c_t)$ and $(F_t(\phi_2))$ are bounded (by a -- non-random
-- constant), and the jumps of $X$ are also bounded by a constant. \qed
\vst

Next, we introduce a number of notation, for which we assume (H)
and heavily use $\si$, as in (\ref{AS2}). Recall $X'=X-X_0-X^c$:
\bee\label{TR1}
\left.\begin{array}{l}
\chi'^n_i=\frac1{\rn}\itai(\si_s-\si_{(i-1)\De_n})~dW_s\\[2mm]
\be^n_i=\si_{(i-1)\de_n}\dd W/\rn,\qquad
\chi^n_i=\chi'^n_i+\frac1{\rn}~\dd X'\\[2mm]
\rho^n_i=\rho_{\si_{(i-1)\De_n}}.
\end{array}\right\}
\eee
In particular, $\dd X=\chi^n_i+\be^n_i$. It is obvious that (SH)
implies for all $q>0$:
\bee\label{TR2}
\left.\begin{array}{l}
\ec(|\be^n_i|^q)\leq K_q,\qquad\ec(|\chi'^n_i|^q)\leq K_q,
\qquad\ec(|\dd X^c|^q)\leq K_q\De_n^{q/2}\\[2mm]
\ec(|\dd X'|^q)\leq \left\{\begin{array}{ll}
K_q\De_n^{1\bigwedge(q/2)}&\mbox{in general}\\
K_q\De_n^q&\mbox{if $X$ is continuous}\end{array}\right.\\[2mm]
\ec(|\chi^n_i|^q)\leq \left\{\begin{array}{ll}
K_q\De_n^{-(1-q/2)^-}&\mbox{in general}\\
K_q&\mbox{if $X$ is continuous}\end{array}\right.\end{array}\right\}
\eee

\begin{lem}\label{LT1} Under (SH) we have
\bee\label{TR3}
\De_n\st\E(\phi_2(\chi^n_i))~\toucp~0.
\eee
\end{lem}

\nib Proof. \rm For any $\ep\in(0,1]$ we write
$X'=N(\ep)+M(\ep)+B(\ep)$, where
$$N(\ep)=(x1_{\{|x|>\ep\}})\star\mu,\quad
M(\ep)=(x1_{\{|x|\leq\ep\}})\is,\quad
B(\ep)=B-(\ka(x)1_{\{|x|>\ep\}})\star\nu.$$
We also set
$$\ga^n_i(y)=\frac1{\De_n}~\ec\left(\itai dt\int_{\{|x|\leq y\}}
\phi_2(x)F_t(dx)\right),$$
which is increasing in $y$ with $\ga^n_i(y)\leq K$ by (SH). Then
$$\pc(\dd N(\ep)\neq0)\leq K\ep^{-2}\De_n,\quad
\ec((\dd M(\ep))^2)\leq\De_n\ga^n_i(\ep),\quad
|\dd B(\ep)|\leq K\De_n\ep^{-1}$$
(use Tchebycheff inequality for the first and last estimates).
We also have
$$\ec((\chi'^n_i)^2)~=~\ga'^n_i~:=~ \frac1{\De_n}~\ec\left(
\itai (\si_u-\si_{(i-1)\De_n})^2~du\right).$$

The following is obvious:
$$\phi_2(\chi^n_i)\leq 1_{\{\dd N(\ep)\neq0\}}+
3|\chi'^n_i|^2+3\De_n^{-1}(|\dd M(\ep)|^2+3\De_n^{-1}|\dd B(\ep)|^2),$$
Then if we take $\ep=\ep_n=\De_n^{1/4}$ we deduce from
the previous estimates that
\bee\label{TR4}
\ec\left(\phi_2(\chi^n_i)\right)
\leq K\rn+K\ga'^n_i +K\ga^n_i(\ep_n) .
\eee

Now, observe that
$$\De_n~\E\left(\st(\ga^n_i(\ep_n)+\ga'^n_i))\right)\leq
\E\left(\int_0^tdu\left((\si_u-\si_{\De_n[u/\De_n]})^2
+\int_{\{|x|\leq\ep_n\}}\phi_2(x)F_u(x)
\right)\right).$$
(SH) implies that for each $(\om,u)$ the middle parenthesis
in the right side above goes to $0$, while staying bounded
by a constant, so by Lebesgue's theorem the left side
goes to $0$. Plugging this into (\ref{TR4}) immediately
gives (\ref{TR3}). \qed

\begin{lem}\label{LT2} Under (SH) we have for all $f\in\ea''_1$
  and all $\rho>0$:
\bee\label{TR5}
\lim_{\ep\to0}~\lim_{A\to\infty}~\limsup_n~\PP\left(\st\ec\Big(
(f(\psi_\ep-\psi_{A\rn}))^2(\dd X)\Big)>\rho\right)=0.
\eee
\end{lem}

\nib Proof. \rm We have $|f(x)|\leq K|x|$ for $|x|\leq 1$, so as soon
as $A\rn\leq \ep/2\leq1/4$ we have (by singling out the two cases
$|x|\leq|y|$ and $|x|>|y|$):
\bean
|f(\psi_\ep-\psi_{A\rn})|(x+y)&\leq&
K|x|1_{\{A\rn/2\leq|x|\leq 3\ep\}}+K|y|1_{\{A\rn/2\leq|y|\leq 3\ep\}}\\
&\leq&
K|x|(\psi_{3\ep}-\psi_{A\rn/2})(x)+K|y|(\psi_{3\ep}-\psi_{A\rn/2})(y).
\eean
Hence it is enough to prove (\ref{TR5}) for $f=h_1$, and
separately for $X^c$ and for $X'$. First, by (\ref{TR2}) we have
$$\ec(|\dd X^c|^2(\psi_\ep-\psi_{A\rn})^2(\dd X^c))\leq
\ec(|\dd X^c|^21_{\{|\dd X^c|\geq A\rn\}})\leq\frac {K\De_n}A,$$
and (\ref{TR5}) for $h_1$ is then obvious for $X^c$. Second, we have
$$\st\ec\Big(|\dd X'|^2(\psi_\ep-\psi_{A\rn})^2(\dd X')\Big)
\leq \st\ec\Big((h_2\psi_\ep)(\dd X'')\Big).$$
Lemma \ref{Lem1} applied to $X=X'$ (note that
$f\psi_\ep$ is bounded) yields that the right side above converges
u.c.p. to $(g\psi_\ep)\star\nu_t$, and the later goes to $0$ as
$\ep\to0$: this shows (\ref{TR5}) for $X'$. \qed

\begin{lem}\label{LT3} Under (H) we have $\De_n\st
\rho^n_{i}(g)\toucp\int_0^t\rho_{\si_s}(g)ds$ if $g\in\ea$ is continuous.
\end{lem}

(The continuity of $g$ is much too strong for this, but the above
result is enough for us). 
\vsq

\nib Proof. \rm The process $\si$ is c\`adl\`ag, hence the function
$s\mapsto\rho_s(g)=\E(g(\si_sU))$ is also c\`adl\`ag by Lebesgue's theorem.
The result is then obvious by Riemann approximation of the integral.  \qed

\begin{lem}\label{LT4} (i) Under (H) any even function $g$ in
  $\ea$ we have
\bee\label{TR6}
\ec\left(\dd N~g(\be^n_i)\right)~=~0
\eee
for $N=W$ and for all $N$ in the set $\n$ of all continuous bounded
martingales which are orthogonal to $W$.

(ii) Assume (SH), and let $g\in\ea$ be
  continuous. If further $q>0$, and $g(x)/|x|^{2/q}\to0$ as
$|x|\to\infty$ when $X$ is not continuous, then
\bee\label{TR7}
\De_n\st\ec\left(\left|g(\dd X/\rn)-g(\be^n_i)\right|^q\right)
~\toucp~0.
\eee
In particular, provided $g(x)/x^2\to0$ as
$|x|\to\infty$ whenever  $X$ is not continuous, then
\bee\label{TR8}
\De_n\BK^n(g)_t~\toucp~\int_0^t\rho_{\si_s}(g)ds.
\eee
\end{lem}

\nib Proof. \rm (i) When $N=W$, we have $\dd Ng(\be^n_i)=
h(\si_{(i-1)\De_n},\dd W)$ for a function $h(x,y)$ which is odd and
with polynomial growth in $y$
when $g$ is even, so obviously (\ref{TR6}) holds. When $N\in\n$,
(\ref{TR6}) is proved in Proposition 4.1 of \cite{BGJPS}.

(ii) Since $\ec(g(\be^n_i))=\rho_i^n(g)$, (\ref{TR8})
readily follows from (\ref{TR7}) for $q=1$ and from Lemma \ref{LT3}.
As for (\ref{TR7}), it amounts to the AN property of the array $(\ze'^n_i)$
defined as follows:
$$\ze'^n_i=\De_n\ec(|\ze^n_i|^q),\qquad
\ze^n_i=g(\dd X/\sqrt{\De_n})-g(\be^n_i).$$

We first prove this result when $g(x)/|x|^{2/q}\to0$ at infinity.
Set $G_A(\ep)=\sup(|g(x+y)-g(x)|:~|x|\leq A,~|y|\leq\ep)$ and
$H_A=\sup_{|x|>A}|g(x)|/|x|^{2/q}$ and $L_A=\sup_{|x|\leq A}|g(x)|$.
We have $G_A(\ep)\to0$ as $\ep\to0$ for all $A$, and $H_A\to0$ as
$A\to\infty$, and $L_A<\infty$ for all $A$, and also $|g(x+y)|\leq
L_B+K_qH_B(|x|^{2/q}+|y|^{2/q})$ for all $B>0$. Then, since
$\dd X/\rn=\be^n_i+\chi^n_i$, and with the notation $W^n_i=\dd W/\rn$,
we obtain for $\ep\in(0,1]$, $A,B>0$, and if $|\si|\leq \Ga$:
$$|\ze^n_i|\leq K
\left(G_A(\ep)+H_B\Big(|\chi^n_i|^{2/q}+|\Ga W^n_i|^{2/q}\Big)+L_B\left(
1_{\{|W^n_i|>A/\Ga\}}+\ep^{-2/q}~\phi_{2/q}(\chi^n_i)\right)\right),$$
and thus by using (\ref{TR2}),
$$\ze'^n_i\leq K\De_n\left(G_A(\ep)^q+H_B^q+L_B^q\PP(|U|>A/\Ga)
+L_B^q\ep^{-2}\ec(\phi_2(\chi^n_i))\right).$$
Therefore if we use (\ref{TR3}) we get
\bean
\st\ze'^n_i&\!\!\!\leq&\!\!\!
Kt\left(G_A(\ep)^q+H_B^q+L^q_B\PP(|U|>A/\Ga)\right)
+KL^q_B\De_n\ep^{-2}\st\ec(\phi_2(\chi^n_i))\\
&&\qquad\qquad~\toucp~Kt\Big(G_A(\ep)^q+H_B^q+L^q_B\PP(|U|>A/\Ga)\Big).
\eean
Then we take $B$ such that $H_B$ is small, then $A$ such that
$L_B\PP(|U|>A/\Ga)$ is small, than $\ep$ such that $G_A(\ep)$ is
small, and we deduce that the array $(\ze'^n_i)$ is AN.

Finally when $X$ is continuous and $g$ is of polynomial growth, we
have $H_A=\infty$, but since now $\chi^n_i=\chi'^n_i$ we can use the
estimate $|g(x)|\leq K(1+|x|^p)$ for some $p>0$ to get
$$|\ze^n_i|\leq K\left(G_A(\ep)+\Big(1+|\chi'^n_i|^p+|W^n_i|^p\Big)\Big(
1_{\{|W^n_i|>A/\Ga\}}+\ep^{-1}~|\chi'^n_i|\Big)\right),$$
hence by H\"older and (\ref{TR2}) we deduce
$$\ze'^n_i\leq K\De_n\left(G_A(\ep)^q+\left(\PP(|U|>A/\Ga)\right)^{1/2}+
\E(|U|^p1_{\{|U|>A/\Ga\}})
+\frac{1}{\ep}~\left(\ec(|\chi'^n_i|^2)\right)^{1/2}\right).$$
Then we may conclude as above, using Lemma 7.8 of \cite{BGJPS} instead
of (\ref{TR4}). \qed

\subsection{Proof of Theorem \ref{T3}.}

This theorem is a consequence of the following two lemmas:
the first one proves the result under the stronger assumptions (SH),
the second one is a standard localization procedure giving the results
under (H).

\begin{lem}\label{LN1} Theorem \ref{T3} holds under the assumption (SH).
\end{lem}

\nib Proof. \rm (i) If $g\in\ea$ is continuous, the process
$$\BV^n(g)_t=\De_n\st\Big(g(\be^n_i)-\rho^n_i(g)\Big).$$
is a square-integrable martingale w.r.t.\ the
filtration $(\f'^n_t=\f_{\De_n[t/\De_n]})_{t\geq0}$, and its
predictable bracket is $\De_n^2\st\Big(\rho^n_i(g^2)-\rho^n_i(g)^2\Big)
\leq Kt\De_n$. Hence $\BV^n(g)\toucp0$ and we deduce from Lemma
\ref{LT3} that
\bee\label{TR16}
\De_n\st g(\be^n_i)~\toucp \int_0^t\rho_u(g)du.
\eee
If further $g(x)/x^2\to0$ as $|x|\to\infty$, or if $X$ is continuous, we
can apply (\ref{TR7}) with $q=1$ to 
deduce via Lenglart's inequality that
$$\De_n\st\Big(g(\dd X/\rn)-g(\be^n_i)\Big)~\toucp~0.$$
Combining this with (\ref{TR16}) gives (i). 

(ii) We have $\De_n^{-r/2}V^n(h_r)=V'^n(h_r)$. If $f\in\ea_r$, for all
$\ep>0$ we also have $|f-h_r|\leq\ep h_{qr}+K_\ep(1-\psi_\ep)h_p$ for
some constants $p>2$ and $K_\ep>0$, hence
\bee\label{TR15}
|\De_n^{1-r/2}V^n(f)-\De_nV'^n(h_r)|\leq \ep\De_nV'^n(h_r)
+\De_n^{1-r/2}V^n((1-\psi_\ep)h_p).
\eee
On the one hand $(1-\psi_\ep)h_p\in\ea_2'''\cap C^0$, so Theorem
\ref{T1} and $r<2$ yield $\De_n^{1-r/2}V^n((1-\psi_\ep)h_p)\toucp0$.
On the other hand $\De_nV'^(h_r)\toucp m_r\int_0^tc_u^{r/2}\,du$ by
(i). Since $\ep>0$ is arbitrarily small, we deduce (ii) from
(\ref{TR15}).  

(iii) For any $0<\ep<A<\infty$ we have $2A\rn\leq
\al\De_n^{\vpi}\leq \ep$ for all $n$ large enough, and if this holds we
have
\bee\label{TR9}
\De_n V'^n(h_2\psi_A)\leq V''^n(\vpi,\al)\leq V^n(h_2\psi_\ep).
\eee
Since $h_2\psi_A$ is bounded continuous, (i) implies that the left
side above converges u.c.p.\ to
$V'(A)_t=\int_0^t\rho_{\si_u}(h_2\psi_A)\,du$, which in turn increases
(u.c.p. again) to $C_t$ as $A\to\infty$. The right side of (\ref{TR9})
has jumps smaller than $4\ep^2$, and since $h_2\psi_\ep\in\ea_2\cap C^0$
it converges in probability for the Skrokhod topology to
$C+(h_2\psi_\ep)\star\mu$ by Theorem \ref{T1}-(b), whereas 
$C+(h_2\psi_\ep)\star\mu$ decreases u.c.p.\ to $C$ as $\ep\to0$. Then
(iii) is obvious. \qed

\begin{lem}\label{LN2} If Theorem \ref{T3} holds under the
  assumption (SH), it also holds under the assumption (H).
\end{lem}

\nib Proof. \rm (H) implies the existence of a sequence of
stopping times $T_p$ increasing to $\infty$ and such that the three
processes $(b_t)$, $(c_t)$ and $(F_t(\phi_2))$ are bounded by a
constant $K_p$ for all $t\leq T_p$, and also such that $|\De X_s|\leq
p$ for all $s<T_p$ (note that we usually cannot find $T_p$ as above,
such that $|\De X_{T_p}|\leq p$). Then the process
$$X(p)_t=X_0+B_{t\bigwedge T_p}+X^c_{t\bigwedge T_p}
+\ka\is_{t\bigwedge T_p}
+(\ka'(x)1_{\{|x|\leq p\}})\star\mu_{t\bigwedge T_p}$$
(compare with (\ref{I1})) satisfies (\ref{AS1}) with $b(p)_t=
b_t1_{\{t\leq T_p\}}$ and $c(p)_t=c_t1_{\{t\leq T_p\}}$ and
$F(p)_t(dx)=1_{\{|x|\leq p\}}\cdot F_t(dx)1_{\{t\leq T_p\}}$, and also
$|\De X(p)|\leq p$ by construction: hence $X(p)$ satisfies (SH).

By hypothesis, for each $p$ the processes $\De_nV'^n(X(p);g)$ in (i)
converge u.c.p. to 
$\int_0^t\rho_{\si(p)_u}(g)du= \int_0^{t\bigwedge
T_p}\rho_{\si_u}(g)du$.  Since we have $V'^n(X(p;g)_t=V'^n(X;g)_t$
for $t<T_p$, whereas $T_p\uparrow\infty$, we readily deduce the result
for (i). For (ii) and (iii) it is proved in the same way. \qed

\section{Proofs for the CLTs}\label{secCLT}
\setcounter{equation}{0}
\renewcommand{\theequation}{\thesection.\arabic{equation}}

\subsection{Technical consequences of (K), (L-$s$) and (H').}

Exactly as for Assumption (H) which was strengthened into (SH) for
technical reasons, we need to strengthen (K), (L-$s$) and (H') as
follows : \vst

\nib Hypothesis (SK): \rm We have (K) and (SH), and
the functions $\ga_k=\ga$ do not depend on $k$ and are bounded. \qed
\vst

\nib Hypothesis (SL-$s$): \rm We have (L-$s$) and the processes
$(b_t)$, $(c_t)$, $(\Wb_t)$, $(\Wsi_t)$, $(\Wsi'_t)$ are bounded, and
the functions $\ga_k=\ga$ and $\Wga_k=\Wga$ do not depend on $k$ and
are bounded. \qed
\vst

\nib Hypothesis (SH'): \rm We have (H) and the process
$(c_t)$ is bounded away from $0$. \qed
\vst

Now we proceed to derive some consequences of these assumptions,
except that the first result, used for Theorem \ref{T4}, needs (SH)
only.

\begin{lem}\label{LTT1} Assume (SH). If $f\in\ea'_1$ we have for all
  $\rho>0$:
\bee\label{TT1}
\lim_{\ep\to0}~\limsup_n~\PP\left(\st\ec
\left(\left|f(\dd X)\psi_\ep(\dd X)-|\rn~\be^n_i|
\right|^2\right)>\rho\right)=0.
\eee
\end{lem}

\nib Proof. \rm Suppose first that $f=h_1$. Observe that for any $A>0$,
$$\left|f(\dd X)\psi_\ep(\dd X)-|\rn~\be^n_i|\right|
\leq \sum_{j)1}^3\ze^n_i(A,\ep,j),$$
where, with the notation $g_A=f\psi_A$,
$$\ze^n_i(A,\ep,1)=\left|f(\dd X)(\psi_\ep(\dd X)
-\psi_{A\rn})(\dd X)\right|,$$
$$\ze^n_i(A,\ep,2)=\rn~\left|g_A(\dd X/\rn)-g_A(\be^n_i)\right|,$$
$$\ze^n_i(A,\ep,3)=\ze^n_i(A,3)=\rn~|\be^n_i|~(1-\psi_A(\be^n_i)).$$
Then it is enough to show that
\bee\label{TT2}
\lim_{\ep\to0}~\lim_{A\to\infty}~\limsup_n~\PP\left(\st\ec
(\ze^n_i(A,\ep,j)^2)>\rho\right)=0
\eee
for all $\rho>0$ and $j=1,2,3$.
Now, (\ref{TT2}) for $j=1$ is exactly (\ref{TR5}), and (\ref{TT2})
for $j=2$ follows from (\ref{TR7}) applied to $g_A$ (which is bounded
continuous) and $q=2$, and (\ref{TT2}) for $j=3$ immediately follows from
(\ref{TR2}).

Finally when $f\in\ea'_1$, in order to get the result it
suffices to prove that the array
$$\ze^n_i(\ep)=\ec\left(\left(~\Big|(f-h_1)(\dd X)
\Big|\psi_\ep(\dd X)\right)^2\right)$$
is AN, for each $\ep>0$. We have $|(f-h_1)\psi_\ep|\leq \eta_\ep
\phi_1$ with $\eta_\ep\to0$ as $\ep\to0$, hence
$$\E\left(\st\ze^n_i\right)\leq \eta_\ep^2~\BH^n_t(\phi_2),$$
and the result follows from Lemma \ref{Lem1}. \qed
\vsc

Note that under (K), resp. (SL-$2$), Equations (\ref{AS3}) and (\ref{AS4})
take the form
\bee\label{TT3}
X_t=X_0+\int_0^tb'_sds+\int_0^t\si_sdW_s+\de\star(\umu-\unu)_t,
\eee
\bee\label{TT4}
\si_t=\si_0+\int_0^t\Wb'_sds+\int_0^t\Wsi_sdW_s+\int_0^t\Wsi'_sdW'_s
+\Wde\star(\umu-\unu)_t,
\eee
where $b'_t=b_t+\int \ka'(\de(t,x))~dx$ and
$\Wb'_t=\Wb_t+\int \ka'(\Wde(t,x))~dx$ are bounded.

The key result is that, under appropriate assumptions, the
convergence in (\ref{TR8}) holds with a rate $1/\rn$. This has
been shown in \cite{BGJPS} when $X$ is continuous, i.e. $\de=0$,
but in general we need some estimate on the increments of the
process $\de\star(\umu-\unu)$. However we start with a result
proved in Section 8-2 of \cite{BGJPS}, for which the process $\de$
plays no role:

\begin{lem}\label{LTT2} Assume (SL-$2$) and let either $g$ be
  differentiable with $g'\in\ea$, or $g\in\ea$ and (SH') hold. Then
$\frac1{\rn}\left(\st\De_n\rho^n_{i}(g) 
-\int_0^t\rho_{\si_s}(f)ds\right)\toucp0$.
\end{lem}

Now we fix a sequence $\ep_n$ in $(0,1)$, going to $0$ and to be
chosen later, and we put $E_n=\{x:\ga(x)>\ep_n\}$. Recall that
$t\mapsto \de(\om,t,x)$ is left continuous with right limits, and
we denote by $\de_+(\om,t,x)$ the right limit at time $t$. Then we
set
\bee\label{TT5} \left.\begin{array}{l}
\ze^n_i(1)=\frac1{\rn}\itai \int_{E_n^c}\de(v,x)(\umu-\unu)(dv,dx)\\
\ze^n_i(2)=-\frac1{\rn}\itai \int_{E_n}\de_+((i-1)\De_n,x)\unu(dv,dx)\\
\ze^n_i(3)=-\frac1{\rn}\itai \int_{E_n}(\de(v,x)-\de_+((i-1)\De_n,x))
\unu(dv,dx)\\
\ze^n_i(4)=\frac1{\rn}\itai \int_{E_n}\de(v,x)\umu(dv,dx),
\end{array}\right\}
\eee
We denote by $\f'^n_i$ the $\si$-field generated by
$\f_{(i-1)\De_n}$ and the variables $(W_u:0\leq u\leq i\De_n)$.

\begin{lem}\label{LTT3} Assume (SL-$s$). The function
$\ga_s(y)=\int_{\{x:\ga(x)\leq y\}}\ga(x)^sdx$ is bounded increasing
and goes to $0$ as $y\to0$, and we have for $r\in(0,1]$ and $\al>0$:
\bee\label{TT6}
\E\left(\ze^n_i(1)^2\mid\f'^n_i\right)\leq K\ep_n^{2-s}\ga_s(\ep_n),
\eee
\bee\label{TT7}
|\ze^n_i(2)|+|\ze^n_i(3)|\leq K\rn~\ep_n^{-(s-1)^+},
\eee
\bee\label{TT8}
\E\left(|\ze^n_i(4)|^r\mid\f'^n_i\right)
\leq K\De_n^{1-r/2}\ep_n^{-(s-r)^+},
\eee
\bee \label{TT9}
\E\left(|\ze^n_i(4)|\bigwedge\al\mid\f'^n_i\right)\leq K\al\De_n\ep_n^{-s}.
\eee
\end{lem}

\nib Proof. \rm Obviously $|\ze^n_i(2)|+|\ze^n_i(3)|\leq 3\rn\int_{E_n}
\ga(x)dx$, hence Tchebycheff's inequality yields (\ref{TT7}).
Next, conditionally on
$\f'^n_i$, the measure $\umu$ restricted to $((i-1)\De_n,\infty)\times\R$ is
still a Poisson measure with intensity measure $\unu$, because $\umu$
and $W$ are independent. Hence
$$\E(\ze^n_i(1)^2\mid\f'^n_i)
=\frac1{\De_n}\E\left(\itai du\int_{E_n^c}\de(u,x)^2dx\Big|
\f'^n_i\right)\leq  \int_{E_n^c}\ga(x)^2dx,$$
and (\ref{TT6}) follows.

Finally $|\ze^n_i(4)|\leq Z^n_i:=\frac1{\rn}\itai\int_{E_n}\ga(x)
\umu(dv,dx)$, and $Z^n_i$ is independent of $\f'^n_i$ and
is a compound Poisson variable. More specifically, if
$\eta_n=\int1_{E_n}(x)dx$, then $Z^n_i$ is the sum of $N$
i.i.d. variables $Y_j$ with
$\E(f(Y_j))=\frac1{\eta_n}\int_{E_n}f(\ga(x)/\rn)dx$ and $N$ is
independent of the $Y_j$'s and Poisson with parameter $\eta_n\De_n$.
We deduce first that
$$\E(|Z^n_i|^r)\leq \E\Big(\sum_{j=1}^N|Y_j|^r\Big)
=\E(N)~\E(|Y_1|^r)=\De_n^{1-r/2}\int_{E_n}\ga(x)^rdx,$$
and second that
$$\E(|Z^n_i|\bigwedge\al)\leq \al\PP(N\geq1)\leq \De_n\eta_n.$$
Since $\int\ga(x)^sdx<\infty$, we deduce (\ref{TT8}) and (\ref{TT9}) from
Tchebycheff's inequality again.  \qed

\begin{lem}\label{LTT4} Under (SL-$s$) we have
\bee\label{TT10}
\rn\st\sqrt{\ec(|\ze^n_i(3)|^2)}~=~\tou(\ep_n^{-(s-1)^+}).
\eee
\end{lem}

\nib Proof. \rm By a repeated application of Cauchy-Schwarz, the
expected value of the left side of (\ref{TT10}), say
$a_n(t)$, satisfies
\bean
&&a_n(t)^2~\leq~t\E\left(\st\ec(|\ze^n_i(3)|^2)\right)\\
&&\leq t\E\left(\st\frac1{\De_n}\left(\itai dv\int_{E_n}
|\de(v,x)-\de_+((i-1)\De_n,x))|x\right)^2\right)\\
&&\leq t\E\left(\int_0^t\!\! dv \int_{E_n}\!
|\de(v,x)-\de_+(\De_n[v/\De_n],x)|^sdx \int_{E_n}\!
|\de(v,x)-\de_+(\De_n[v/\De_n],x)|^{2-s}dx\right)\\
&&\leq t\eta_n(s)\E\left(\int_0^t dv\int_{E_n}
|\de(v,x)-\de_+(\De_n[v/\De_n],x)|^sdx\right),
\eean
where $\eta_n(s)=\int_{E_n}(2\ga(x))^{2-s}dx$. Now, on the one hand
the expectation in the last term above goes to $0$ because of the
properties of $\de$ and of (SL-$s$), by an application of Lebesgue's
theorem. On the other hand we have $\eta_n(s)\leq K$ is $s\leq1$, and
when $s>1$ we have $\eta_n(s)\leq K\ep_n^{2-2s}$ by Tchebycheff's
inequality, and the result follows. \qed
\vsc

We are now ready to improve on (\ref{TR8}) by giving a rate, at least
in some special situations. That is, we give estimates on the processes
\bee\label{TR8'}
U^n(g)_t~=~\De_n\BK^n(g)_t-\int_0^t\rho_{\si_u}(g)~du
\eee
in three different situations\,:

$\bullet$ Case (i): $g$ is $C^2_b$,

$\bullet$ Case (ii): $g_n=h_2\,\psi_{\al\De_n^{\vpi-1/2}}$ for some $\al>0$ and
some $\vpi\in(0,1/2)$,

$\bullet$ Case (iii): $g=h_r$ for some $r\in(0,1)$,

\noindent and also when $g$ is $C^1$ with $g'\in\ea$, when $X$ is continuous. 
We also need to introduce the following functions $\eta_s$ on
$[0,1]$, where $s\in[0,2]$:
\bee\label{TT11}
\eta_s(r)=\frac{(2-s)(1+r)(2-r)}{4+2s(1-r)}.
\eee

\begin{lem}\label{LTT5} Assume (SL-$s$). Then 
\bee\label{TT111}
\frac1{\rn}~U^n(g)~\toucp~0,\qquad\mbox{or }~~~
\frac1{\rn}~U^n(g_n)~\toucp~0
\eee
in the following cases:

(a) $X$ is continuous and either $g$ is $C^1$ with $g'\in\ea$ and $g$
is even, or $g=h_r$ for $r\in(0,1]$ if further (SH') holds;

(b) in case (i) if~ $s\leq1$;

(c) in case (ii) if~ $0\leq s\leq \frac{4\vpi-1}{2\vpi}$ (hence
$\vpi\geq\frac14$ and $s<1$); 

(d) in case (iii) if further (SH') holds, and provided either $s\leq
\frac23$ and $0<r<1$,
or $\frac23<s<1$ and $\frac{1-\sqrt{3s^2-8s+5}}{2-s}<r<1$,

\noindent Otherwise, we have for all $\ep>0$:
\bee\label{TT12}
U^n(g_n)_t=
\left\{\begin{array}{ll}
\tou(\De_n^{1-s/2})&\mbox{in case (i) if }~s>1 \\[1.5mm]
\tou(\De_n^{(2-s)\vpi})&\mbox{in case (ii) if }~
s>\frac{4\vpi-1}{2\vpi}\\[1.5mm]
\tOu(\rn)\quad&\mbox{if $~g_n=h_1$, (SH') holds and }~s\leq 1\\[1.5mm]
\tou(\De_n^{\eta_s(r)-\ep})&\mbox{if $~g_n=h_r$ and (SH') holds
  and }\\&\mbox{  and either}~0<r\leq1<s,~~\mbox{or $~\frac23<s<1~$ and}\\
&\mbox{$~0<r\leq \frac{1-\sqrt{3s^2-8s+5}}{2-s}~$ and $~r<1$.}
\end{array}\right.
\eee
\end{lem}

In the last case of (\ref{TT12}) we have $0<\eta_s(r)<1/2$.
Comparing with (5.14) of \cite{BGJPS}, the case (b) above and
the first estimate in (\ref{TT12}) are just as
good, except for the regularity conditions on $g$; but we suspect that
the last two estimates in (\ref{TT12}) are not optimal: when $X$ is
a L\'evy process we have $U^n(h_r)_t=\tou(\rn)$ 
as soon as $s\leq1$ and $r<1$.
The same comment will also apply to Lemma \ref{LTT6} below. Recall that
when $X$ is continuous the assumptions (SL-$s$) for $s\in[0,2]$
are all equivalent.
\vsq

\nib Proof. \rm Throughout, we assume (SL-$s$). Since (a) is
in \cite{BGJPS}, we only consider (b), (c), (d) and (\ref{TT12}), with
$g_n=g$ in cases (i) and (ii), and we set $U^n_t=U^n(g_n)_t$ and
$\al_n=\al\De_n^{\vpi-1/2}$. The proof goes through several steps.

a) First, we state some obvious properties of our functions $g_n$\,:
\bee\label{TT13}
\left.\begin{array}{l}
|g(x+y)-g(x)|\leq K(|y|\bigwedge1)\\
|g(x+y)-g(x)-g'(x)y|\leq K(|y|\bigwedge y^2)
\end{array}\right\}\quad\mbox{in case (i)},
\eee
\bee\label{TT13'}
\left.\begin{array}{l}
|g_n(x+y)-g_n(x)|\leq K\al_n(|y|\bigwedge\al_n)\\
|g_n(x+y)-g_n(x)-g'_n(x)y|\leq Ky^2
\end{array}\right\}\quad\mbox{in case (ii)},
\eee
\bee\label{TT14}
\left.\begin{array}{l}
|g_n(x+y)-g_n(x)|\leq K|y|^r\\
x\neq0\quad\Rightarrow\quad |g'_n(x)|\leq K|x|^{r-1}\\
0<|y|\leq\frac{|x|}2\quad\Rightarrow\quad |g'_n(x+y)-g'_n(x)|\leq K|x|^{r-2}|y|
\end{array}\right\}\quad\mbox{in case (ii)},
\eee

b) Recall (\ref{TR1}) and set $\be'^n_i=\be^n_i+\chi^n_i-\ze^n_i(4)$. Using
the previous estimates, we readily deduce from Lemma \ref{LTT3} that
(recall the notation (\ref{TT5})):
\bee\label{TT15}
  \st\ec\left(\Big|g_n(\be^n_i+\chi^n_i)-g_n(\be'^n_i)\Big|\right)
=\left\{\begin{array}{ll}
\tOu(\ep_n^{-s})&\mbox{in case (i)}\\[1.5mm]
\tOu(\al_n^2\ep_n^{-s})&\mbox{in case (ii)}\\[1mm]
 \tOu(\De_n^{-r/2}\ep_n^{-(s-r)^+})&\mbox{in case (iii)}
\end{array}\right.
\eee

c) Now we state some results of \cite{BGJPS}, see in particular
Lemma 7.7\,: put $\wbe^n_i=1+|\be^n_i|^{-1}$ and $Z^n_i=1+|\be^n_i|^q$
for some $q\geq0$. Then under (SH'), for any $\te\in(1,2)$ and
$l\in(0,1)$ the variables $\chi^n_i-\sum_{j=1}^4\ze^n_i(j)$ (in which
the jumps of $X$ play no role) are of the form $\wxi^n_i+\Wxi^n_i$, a
decomposition which depends on $\te$ and $l$ and which satisfies
\bee\label{TT16}
\left.\begin{array}{l}
\st\ec(Z^n_i~|\wbe^n_i|^l~|\Wxi^n_i|^\te)=\tOu(\De_n^{\te/2-1})\\[2mm]
\st\ec(Z^n_i~|\wbe^n_i|^l~|\wxi^n_i|)=\tou(\De_n^{-1/2}),
\end{array}\right\}
\eee
and also for all odd functions $k$ in $\ea$ :
\bee\label{TT17}
\ec(k(\be^n_i)\Wxi^n_i)=0.
\eee
Furthermore a look at the proof of the afore-mentioned lemma shows
that when $l=0$ these hold also without (SH') and for $\te=2$, and
that (again without (SH'))
\bee\label{TT16'}
\st\ec(|\wxi^n_i|^2)=\tOu(1).
\eee

In cases (i) we set $\te=2$ and $l=0$ and $Z^n_i=1$. In case (ii) we
set $\te=2$ and $l=0$ and $Z^n_i=1+|\be^n_i|$. In case (iii) 
we fix $\te\in(1,r+1)$ (to be chosen
later), and set $l=\te-r$ which is in $[0,1)$, and
$Z^n_i=1+|\be^n_i|^{\te-1}$. 

Now we set $\Wxi'^n_i=\Wxi^n_i+\sum_{j=1}^3\ze^n_i(j)$, so that
$\be'^n_i=\be^n_i+\wxi^n_i+\Wxi'^n_i$. Since $\ec(Z^n_i~|\wbe^n_i|^l)
\leq K$ when $l$ is as above, and since $Z^n_i$ and $\wbe^n_i$ are
$\f'^n_i$-measurable, we deduce from (\ref{TT6}), (\ref{TT7}) and
(\ref{TT16}) that (with $l=0$ in cases (i,ii), and $l$ as above, under
(SH'), in case (iii))\,:
\bee\label{TT18}
\st\ec(Z^n_i~|\wbe^n_i|^l~|\Wxi'^n_i|^\te)=\tOu\left(\De_n^{\te/2-1}
\ep_n^{-\te(s-1)^+}+\De_n^{-1}\ep_n^{\te(2-s)/2}\ga_s(\ep_n)^{\te/2}\right).
\eee
If $k\in\ea$ is odd and since $\ze^n_i(2)$ is $\f_{(i-1)\De_n}$-measurable,
we clearly have $\ec(g'(\be^n_i)\ze^n_i(j))=0$ for $j=2$, whereas this also
holds for $j=1$ because, as in the proof of Lemma \ref{LTT3}, we have
$\E(\ze^n_i(1)\mid\f'^n_i)=0$. Furthermore we have $|\ec(k(\be^n_i)
\ze^n_i(3))|\leq K\sqrt{\ec(|\ze^n_i(3)|^2)}$ by Cauchy-Schwarz, and
because $\ec(k(\be^n_i)^2)\leq K$ by (\ref{TR2}). Combining these facts
with (\ref{TT17}) and (\ref{TT10}), we obtain
\bee\label{TT19}
\st\left|\ec(k(\be^n_i)\Wxi'^n_i)\right|~=~\tou(\De_n^{-1/2}\ep_n^{-(s-1)^2}).
\eee

d) Now, to evaluate $g_n(\be'^n_i)-g_n(\be^n_i)$, we partly reproduce a
proof in \cite{BGJPS}, with some relevant changes. We first treat case
(iii), which is the most difficult. Set
$A^n_i=\{|\Wxi'^n_i+\wxi^n_i|>|\be^n_i|/2\}$. By a Taylor expansion, we have
$g(\be'^n_i)-g(\be^n_i)=g'(\be^n_i)\Wxi'^n_i+\sum_{j=1}^4\de^n_i(j)$, where
$$\de^n_i(1)=(g(\be'^n_i)-g(\be^n_i))1_{A^n_i},\qquad
\de^n_i(2)=-g'(\be^n_i)(\be'^n_i-\be^n_i)1_{A^n_i},$$
$$\de^n_i(3)=g'(\be^n_i)\wxi^n_i,,\qquad
\de^n_i(4)=(g'(\be''^n_i)-g'(\be^n_i))(\be'^n_i-\be^n_i)1_{(A^n_i)^c},$$
where $\be''^n_i$ is a random variable which is between $\be^n_i$ and
$\be'^n_i$.

If we apply (\ref{TT14}) and single out the two
cases $|\Wxi'^n_i|\geq|\wxi^n_i|$ and $|\Wxi'^n_i|<|\wxi^n_i|$, upon
observing that in the first case for instance we have
$|\Wxi'^n_i|>|\be^n_i|/4$ if we are on $A^n_i$, we get
\bean
|\de^n_i(1)|&\leq& K|\Wxi'^n_i+\wxi^n_i|^r)1_{A^n_i}
\\&\leq& K|\wxi^n_i|~|\be^n_i|^{r-1}+K|\Wxi'^n_i|^\te~
|\be^n_i|^{r-\te}~\leq~
KZ^n_i~|\wbe^n_i|^l\left(|\wxi^n_i|+|\Wxi'^n_i|^{\te}\right)
\eean
(recall $\te\geq r$). In a similar way,
$$|\de^n_i(2)|+|\de^n_i(3)|~\leq~
K|\wxi^n_i|~|\be^n_i|^{r-1}+K|\Wxi^n_i|^{\te}~|\be^n_i|^{r-\te}
~\leq~KZ^n_i~|\wbe^n_i|^l\left(|\wxi^n_i|+|\Wxi'^n_i|^{\te}\right).$$
Finally, by singling out the
cases $|\Wxi'^n_i|\geq|\wxi^n_i|$ and $|\Wxi'^n_i|<|\wxi^n_i|$ once more,
\bean
|\de^n_i(4)|&\leq& K|\be^n_i|^{r-2}~
(\Wxi'^n_i+\wxi^n_i)^21_{\{|\Wxi'^n_i+\wxi^n_i|\leq|\be^n_i|/2\}}\\
&\leq& K|\be^n_i|^{r-2}~\left(|\wxi^n_i|~|\be^n_i|
+|\Wxi'^n_i|^{\te}~|\be^n_i|^{2-\te}\right)~\leq~
KZ^n_i~|\wbe^n_i|^l\left(|\wxi^n_i|+|\Wxi'^n_i|^{\te}\right).
\eean
Put these three estimates together to get
\bee\label{TT20}
\Big|g(\be'^n_i)-g(\be^n_i)-g'(\be^n_i)\Wxi'^n_i\Big|\leq KZ^n_i
~|\wbe^n_i|^{l}~\Big(|\wxi^n_i|+|\Wxi'^n_i|^\te~\Big).
\eee

In case (i) things are easier. Indeed, (\ref{TT13}) implies:
\bea
\Big|g(\be'^n_i)-g(\be^n_i)-g'(\be^n_i)\Wxi'^n_i\Big|
&\leq& K(|\Wxi'^n_i+\wxi^n_i|^2\bigwedge|\Wxi'^n_i+\wxi^n_i|)+
K|\wxi^n_i|\nonumber\\
&\leq& K(|\wxi^n_i|+|\Wxi'^n_i|^2).\label{TT21}
\eea

In case (ii), we use the fact that $|g'_n(x)|\leq K|x|)$ and
(\ref{TT13'}) to get :
\bea
\Big|g_n(\be'^n_i)-g_n(\be^n_i)-g'_n(\be^n_i)\Wxi'^n_i\Big|
&\leq& K|\Wxi'^n_i+\wxi^n_i|^2+K|\be^n_i\wxi^n_i|\nonumber\\
&\leq& K(|\wxi^n_i|^2+|\Wxi'^n_i|^2+Z^n_i|\wxi^n_i|).\label{TT21'}
\eea

Moreover in all cases $g_n$ is even, hence $g'_n$ is odd. Therefore
(\ref{TT16}), (\ref{TT16'}), (\ref{TT18}) and (\ref{TT19}), together
with either (\ref{TT20}) or (\ref{TT21}), imply that in all cases
(recall $\te=2$ in cases (i) and (ii)):
\bean
&&\De_n\st\Big|\ec(g_n(\be'^n_i)-g_n(\be^n_i))\Big|\\
&&\qquad\qquad=
\tou\left(\De_n^{1/2}\ep_n^{-(s-1)^+}\right)
+\tOu\left(\De_n^{\te/2}\ep_n^{-\te(s-1)^+}
+\ep_n^{\te(2-s)/2}\ga_s(\ep_n)^{\te/2}\right).
\eean

d) The previous estimate plus Lemma \ref{LTT2} and (\ref{TT15}) yield
$$U^n_t=\left\{\begin{array}{ll}
\tOu(\De_n\ep_n^{-2(s-1)^+}+\De_n\ep_n^{-s}+\ep_n^{2-s}\ga_s(\ep_n))&\\[1.5mm]
\qquad+\tou(\rn~\ep_n^{-(s-1)^+})\hskip3.5cm&\mbox{in case (i)}\\[1.5mm]
\tOu(\De_n\ep_n^{-2(s-1)^+}+\De_n\al_n^2\ep_n^{-s}
+\ep_n^{2-s}\ga_s(\ep_n))&\\[1.5mm]
\qquad+\tou(\rn~\ep_n^{-(s-1)^+})\hskip3.5cm&\mbox{in case (ii)}\\[1.5mm]
\tOu(\De_n^{\te/2}\ep_n^{-\te(s-1)^+}+\De_n^{1-r/2}\ep_n^{-(s-r)^+}
+\ep_n^{\te(2-s)/2}\ga_s(\ep_n)^{\te/2})&\\[1.5mm]
\qquad+\tou(\rn~\ep_n^{-(s-1)^+})\hskip3.5cm&\mbox{in case (iii).}
\end{array}\right.$$
So it remains to choose appropriately $\te\in[1,1+r)$ in case (iii) and
the sequence $\ep_n$ in all three cases.

In case (i) with $s<1$ take $\ep_n=\rn$, to obtain
$Y^n_t=\tou(\rn)$. In case (i) with $s\geq1$ take $\ep_n=A\rn$ for
some $A>1$, to obtain for some $K$ not depending on $A$:
$$\frac{|U^n_t|}{\De_n^{1-s/2}}\leq \left\{\begin{array}{ll}
K\left(\rn+A^{-1}+A^{2-s}\ga_s(A\rn)\right)+\tou(1)\quad
&\mbox{if }~s=1\\
K\left(A^{1-s}+B^{2-s}\ga_s(A\rn)\right)&\mbox{if }~s>1.\end{array}\right.$$
Since $\ga_s(A\rn)\to0$ for any $A$ and since $K$ is arbitrarily
large, we readily deduce that $U^n_t=\tou(\De_n^{1-s/2})$, and we
have the result.

Next, in case (i) we take $\ep_n=A\De_n^\vpi$ for
some $A>1$, to obtain for some $K$ not depending on $A$:
$$|U^n_t|\leq 
K\De_n^{\vpi(2-s)}\left(A^{-s}+A^{2-s}\ga_s(A\De_n^\vpi)\right)
+\tou\left(A^{-(s-1)^+}\De_n^{1/2-\vpi(s-1)^+}\right),$$
and we conclude as in case (i).

Finally consider case (iii). When $s\leq r$ we choose $\te$ arbitrarily in
$(1,1+r)$ and $\ep_n=\De_n^{1/\te(2-s)}$, and the result in (c) is obvious.
When $r<s$, for any given $\te\in(1,1+r)$ we choose
$\ep_n=\De_n^{(2-r)/((2-s)\te+2(s-r))}$. After a simple (although a
bit tedious) computation we deduce
$$U^n_t=\tOu\left(\De_n^{a(\te,r,s)}\right)
+\tou\left(\De_n^{a'(\te,r,s)}\right),$$
with $a(\te,r,s)=\frac{(2-r)(2-s)\te}{2((2-s)\te+2(s-r))}$ and
$a'(\te,r,s)=\frac12$ if $s\leq1$ and
$a'(\te,r,s)=\frac{(2-s)(\te+2(1-r))}{2((2-s)\te+2(s-r))}$ if
$s>1$. Observe that $a'(\te,r,s)\geq a(\te,r,s)$ is $s>1$ and $\te<2$,
and that $a$ is increasing in $\te$. So we should $\te$ as big as
possible (with $\te<1+r$, though), and another simple computation shows
that when $s\leq1$, then $a(1+r,r,s)>1/2$ if and only if
$\frac{1-\sqrt{5-8s+3s^2}}{2-s}<r\leq1$. Then again we get
the results. \qed

\begin{lem}\label{LTT6} Assume (SL-$s$) and let
$g_n$ be as in case (ii) of Lemma \ref{LTT5}. Then for all $t>0$ we have 
\bee\label{TT27}
\st\ec\left(\left(g_n(\dd X/\rn)-(\be^n_i)^2\right)^2\right)
~=~\tou\left(\De_n^{4\vpi-2-s\vpi}\right).
\eee
\end{lem}

\nib Proof. \rm First, we have
\bean
\st\ec\left(\left((\be^n_i)^2-g_n(\be^n_i)\right)^2\right)
&\leq& \st\ec\left((\be^n_i)^41_{\{|\be^n_i|>\al\De_n^{\vpi-1/2}\}}
\right)\\
&\leq&K\De_n^{q((1/2-\vpi)}\st\ec(|\be^n_i|^{4+q})
~\leq~Kt\De_n,
\eean
where the second inequality is valid for all $q>0$ and the third one
is obtained by choosing $q=\frac4{1-2\vpi}$, and using the boundedness of
$\si$. Therefore is it enough to prove 
\bee\label{TT28}
a_n(t)~:=~\st\ec\left(\left(g_n(\dd X/\rn)-g_n(\be^n_i)\right)^2\right)
~=~\tou\left(\De_n^{4\vpi-1-s\vpi}\right).
\eee

With the notation of the proof of Lemma \ref{LTT5} we have
$\dd X=\be^n_i+\be'^n_i+\ze^n_i(4)$, hence by (\ref{TT13'}) we obtain
$$\left(g_n(\dd X/\rn)-g_n(\be^n_i)\right)^2\leq
K\al_n^2|\be'^n_i|^2+K\al_n^3(|\ze^n_i(4)\bigwedge\al_n).$$
Recall also that $\be'^n_i=\wxi^n_i+\Wxi'^n_i$, and we have
(\ref{TT16'}), and also (\ref{TT18} with $\te=2$ and $l=0$ and
$Z^n_i=1$. Therefore we easily deduce from these two estimates, plus
(\ref{TT9}), that (recall $\al_n=\al\De_n^{\vpi-1/2}$)
$$a_n(t)=\tOu\left(\De_n^{4\vpi-2}\ep_n^{-s}+\De_n^{2\vpi-1}\ep_n^{-(s-1)^+}
+\De_n^{2\vpi-2}\ep_n^{2-s}\ga_s(\ep_n)\right).$$
It remains to choose the sequence $\ep_n$, and we take $\ep_n=A\De_n^\vpi$ for
some $A>1$, which gives
$$a_n(t)\leq 
K\De_n^{4\vpi-2-s\vpi}\left(A^{-s}+A^{2-s}\ga_s(A\De_n^\vpi)
+\De_n^{1-2\vpi+\vpi(s\wedge1)}A^{((s-1)^+}\right),$$
and we conclude as in the end of the previous proof. \qed

\subsection{An auxiliary CLT.}

We first give a sketchy proof for a result which is essentially known
already, and which is a CLT for processes of the form
\bee\label{CD1}
\BU^n(g)_t=\rn~\st\Big(g(\be^n_i)-\rho^n_i(g)\Big).
\eee
The assumption (SH) below is of
course much too strong for the result. For the needs of Theorem
\ref{T8} later one, we give a multidimensional version\,: 

\begin{lem}\label{LCONV} Assume (SH) and let $g_1,\cdots,g_d$ be
continuous even functions in $\ea$. The $d$-dimensional
processes $\BU^n$ with components $\BU^n(g_j)$ converge stably
in law to a limit $\BU$ with components $\BU^j_t=\sum_{k=1}^d
\int_0^t\te_u^{jk}~d\BW^k_u$,
where the process $\te$ is $(\f_t)$-optional $d\times d$-dimensional
processes which satisfies
\bee\label{CC3}
(\te_t\te_t^\star)^{jk}~=~
\rho_{\si_t}(g_jg_k)-\rho_{\si_t}(g_j)\rho_{\si_t}(g_k).
\eee
\end{lem}

\nib Proof. \rm Observe that $\BU^n_t=\st\ze^n_i$, where $\ze^n_i$
is the $\f_{i\De_n}$-measurable variable with
components $\ze_i^{n,j}=\rn~\Big(g_j(\be^n_i)-\rho^n_i(g_j)\Big)$.
Moreover $\ec(\ze^n_i)=0$ and $\ec(\|\ze^n_i\|^4)\leq K\De_n^2$
(because $\si_t$ is bounded).
Then a criterion for the stable convergence in law, which can be
found in Theorems IX.7.19 and IX.7.28 of \cite{JS}, gives us the
result, provided we have the following two properties:
\bee\label{CC4}
\st\ec(\ze^{n,j}_i\ze_i^{n,k})~\toucp~\int_0^t(\te_u\te_u^\star)^{jk}~du,
\eee
\bee\label{CC6}
\st\ec(\ze^{n,j}_i~\dd N)\toucp0,\qquad\mbox{if $N=W$ or if }~N\in\n
\eee
(recall that $\n$ is the set of all bounded $(\f_t)$-martingales which
are orthogonal to $W$). (\ref{CC6}) follows from (\ref{TR6}) because
$g_j$ is even. Finally $\ec(\ze^{n,j}_i\ze_i^{n,k}))$ equals the right
side of (\ref{CC3}) evaluated at time $u=(i-1)\De_n$, and multiplied
by $\De_n$. Since the right side of (\ref{CC3}) is a c\`adl\`ag
function of $u$, (\ref{CC4}) follows from Riemann approximation of the
integral on the right, and we are done. \qed
\vsc

For the purpose of proving Theorem \ref{T7}-(ii) we need more than
this lemma. Suppose that we have (SK). Then (\ref{AS3}) holds with
$|\de(\om,t,x)|\leq\ga(x)\leq K$ and $|\Wde(\om,t,x)|\leq\ga(x)\leq K$
and $\int \ga(x)^2dx<\infty$.

We fix $\ep>0$ and consider the process $N=1_E\star\umu$, where
$E=\{x:\ga(x)>\ep\}$. Hence $N$ is a Poisson process with parameter
the Lebesgue measure of $E$, say $\la$. We introduce some notation similar to
(\ref{KR1}), and which depends on $\ep$ :
\bee\label{C4}
\left.\begin{array}{l}
\bullet~ S_1,S_2,\cdots~~\mbox{are the successive jump times of
  $N$},\\
\bullet~ I(n,p)=i,~~~S_-(n,p)=(i-1)\De_n,~~~S_+(n,p)=i\De_n\\
\qquad\qquad\mbox{on the set}~~ \{(i-1)\De_n<S_p\leq i\De_n\},\\
\bullet~ \al_-(n,p)=\frac1{\rn}\Big(W_{S_p}-W_{S_-(n,p)}\Big),
\quad \al_+(n,p)=\frac1{\rn}\Big(W_{S_+(n,p)}-W_{S_p}\Big)\\
\bullet~ R_p=\De X_{S_p},\\
\bullet~ X(\ep)_t=X_t-\sum_{p:~S_p\leq t}R_p,\\
\bullet~ R'^n_p=\dd X(\ep)~~ \mbox{on the set }~
\{(i-1)\De_n<S_p\leq i\De_n\},\\
\bullet~ R'_p=\sqrt{\ka_p}~U_p~\si_{S_p-}+\sqrt{1-\ka_p}~U'_p~\si_{S_p},\\
\bullet~ \Om_n(T,\ep)~~\mbox{is the set of all $\om$ such that each
  interval}~[0,T]\cap((i-1)\De_n,i\De_n]\\
\quad~\mbox{contains at most one $S_p(\om)$, and
that $|\dd X(\ep)(\om)|\leq2\ep$ for all $i\leq T/\De_n$.}
\end{array}\!\!\right\}
\eee
We also suppose that we have the functions $g_j$ and the processes
$\BU^n$ and $\BU$ of Lemma \ref{LCONV}. Then we have the following result,
which is
very close to Lemma 6.2 of \cite{JP}, but with a more involved proof
because we want no restriction on the $\si$-fields $\f_t$.

\begin{lem}\label{LC5} Under (SK), the sequences
  $\left(\BU^n,(\al_-(n,p),\al_+(n,p))_{p\geq1}\right)$
converge stably in law to $\left(\BU,(\sqrt{\ka_p}~U_p~,
\sqrt{1-\ka_p}~U'_p)_{p\geq1}\right)$ as $n\to\infty$.
\end{lem}

\nib Proof. \it Step 1. \rm We need to prove the following : for
all bounded
$\f$-measurable variables $\Psi$ and all bounded Lipschitz functions
$\Phi$ on the Skorokhod space of $d$-dimensional functions on $\R_+$
endowed with a distance for the Skorokhod topology, and all
$q\geq1$ and all continuous bounded functions $f_p$ on $\R^2$, and
with $A_n=\prod_{p=1}^qf_p(\al_-(n,p),\al_+(n,p))$, then
\bee
\E\left(\Psi~\Phi(\BU^n)~A_n\right)
~\to~\WE(\Psi~\Phi(\BU))~\prod_{p=1}^q
\WE\Big(f_p(\sqrt{\ka_p}~U_p~,\sqrt{1-\ka_p}~U'_p)\Big).\label{C40}
\eee
Up to substituting $\Psi$ with $\E(\Psi|\g)$ in both sides, it is
enough to prove this when $\Psi$ is measurable w.r.t.\ the separable
$\si$-field $\g$ generated by the measure $\mu$ and the processes $b$,
$\si$, $W$ and $X$.
\vst

\it Step 2. \rm
We denote by $\umu'$ and $\umu''$ (resp. $\unu'$ and $\unu''$)
respectively the restrictions of $\umu$ (resp. $\unu$) to $\R_+\times
E^c$ and to $\R_+\times E$. We also denote by $(\f'_t)$ the smallest
filtration containing $(\f_t)$ and such that the measure $\umu''$ is
$\f'_0$-measurable. Then $W$ and $\umu'$ are a Wiener process and a
Poisson measure with compensator $\unu'$, relative to $(\f_t)$ of
course, but also to $(\f'_t)$.

Next, for any integer $m\geq1$ we set $S_p^{m-}=(S_p-1/m)^+$ and
$S_p^{m+}=S_p+1/m$, and $B_m=\cup_{p\geq1}(S_p^{m-},S_p^{m+}]$. The
indicator function $1_{B_m}(\om,t)$ is $\f'_0\otimes\rR_+$-measurable,
so the stochastic integral 
$W(m)_t=\int_0^t1_{B_m}(u)~dW_u$ is well defined. We call $(\f'^m_t)$
the smallest filtration containing $(\f'_t)$ and such that the process
$W(m)$ is $\f'^m_0$-measurable, and $\Ga_n(m,t)$ the
set of all integers $i\geq1$ such that $i\leq[t/\De_n]$ and that
$B_m\cap((i-1)\De_n,i\De_n]=\emptyset$, and we introduce the
$d$-dimensional processes $\BU'^n(m)$ and $\BU(m)$ (with $\te$ as in
(\ref{CC3})) with components :
$$\BU'^{n,j}(m)_t~=~\sum_{i\in\Ga_n(m,t)}
\Big(g_j(\be^n_i)-\rho^n_i(g_j)\Big),\quad
\BU^j(m)_t=\sum_{k=1}^d\int_0^t\te^{jk}_u1_{B_m^c}(u)~d\BW^k_u.$$
Again, the integrals above are well defined because $\BW$ is a
Brownian motion w.r.t.\ the smallest filtration containing $(\Wf_t)$
and also $\f'^m_0$ at time $0$. Furthermore $B_m$ decreases to the
union of the graphs of the $S_p$'s, hence $\BU(m)\toucp\BU$ as 
$m\to\infty$. We also have for some $p>0$ because $g_j\in\ea$ and
$\si$ is bounded:
\bean
\E\left(\sup_{s\leq t}|\BU^n(g_j)_s-\BU'^{n,j}(m)_s|^2\right)
&\leq&\E\left(\sum_{p\geq1}\sum_{i:~i\De_n\leq t,|i\De_n-S_p|\leq2/m}
(g_j(\be^n_i)-\rho^n_i(g_j))\right)\\
&\leq& K\E\left(\sum_{p\geq1}\sum_{i:~i\De_n\leq t,|i\De_n-S_p|\leq2/m}
(1+|\dd W|^p)\right)\\
&\leq& \frac Km~\E(\sum_{p=1}^\infty 1_{\{S_p\leq t+1\}})
~\leq~ \frac{Kt}m.
\eean
Therefore, since $\Phi$ is Lipschitz and bounded, it is clearly enough to prove
\bee
\E\left(\Psi~\Phi(\BU'^n(m))~A_n\right)
\to\WE(\Psi~\Phi(\BU(m)))~\prod_{p=1}^q
\WE\Big(f_p(\sqrt{\ka_p}~U_p~,\sqrt{1-\ka_p}~U'_p)\Big)  \label{C41}
\eee
for each $m$, and for $\Psi$ being $\g$-measurable bounded.
\vst

\it Step 3. \rm In the sequel we fix $m$, and we introduce a regular
version $Q=Q_\om(.)$ of the probability $\PP$ on $(\Om,\g)$,
conditional on $\f'^m_0$, and accordingly $\WQ=Q\otimes\PP'$.

Since $\dd W$ is independent of $\f'^m_0$ when $i\in\Ga_n(m,t)$ it is also
standard normal under each $Q_\om$, and in
the proof of Lemma \ref{LCONV} we can replace $\ec$ by the conditional
expectation $\E_{Q_\om}(.|\f'^m_{(i-1)\De_n})$. Moreover $B_m^c$ is a locally
finite union of intervals, hence we still have the convergence in
(\ref{CC4}) if the sum on the left is taken over $\Ga_n(m,t)$ and on
the right we plug in $1_{B_m^c}$ in the integral. Hence
$\BU'^n(m)\tols \BU(m)$ under the measure $Q_\om$, that is
\bee\label{C42}
\E_{Q_\om}(\Psi~\Phi(\BU'^n(m))~\to~\E_{\WQ_\om}(\Psi~\Phi(\BU(m)).
\eee

\it Step 4) \rm By construction $A_n$ is $\f'^m_0$-measurable, so the
left side of (\ref{C41}) is
\bean
&&\E\Big(A_n~\E_{Q_.}(\Psi~\Phi(\BU'^n(m))\Big)=
\E\Big(A_n~\E_{\WQ_.}(\Psi~\Phi(\BU(m))\Big)\\
&&\hskip2cm +\E\Big(A_n~\Big(\E_{Q_.}(\Psi~\Phi(\BU'^n(m))-
\E_{\WQ_.}(\Psi~\Phi(\BU(m))\Big)\Big).
\eean
Since everything above is bounded, the second summand on the right
goes to $0$ by (\ref{C42}), whereas $\Psi'=
\E_{\WQ_.}(\Psi~\Phi(\BU(m))$ is another bounded $\f'^m_0$-measurable
variable. Hence (\ref{C41}) amounts to proving
$$\E\left(\Psi~A_n\right)
~\to~\E(\Psi)~\prod_{p=1}^q
\E\Big(f_p(\sqrt{\ka_p}~U_p~,\sqrt{1-\ka_p}~U'_p)\Big),$$
which is exactly
$(\al_-(n,p),\al_+(n,p))_{p\geq1}\tols(\sqrt{\ka_p}~U_p~,
\sqrt{1-\ka_p}~U'_p)_{p\geq1}$ as $n\to\infty$.
But now, this is a consequence of Lemma 6.2 of \cite{JP} in a slightly simpler
situation, namely we replace $\al^n_j$ and $\be^n_j$ in that lemma by
$\al_-(n,p)$ and $\al_+(n,p)$ here, respectively, and we do not
consider the process $H^{n,\ep}$ in it. Hence we are done. \qed

\begin{lem}\label{LC6} Under the assumptions of Lemma
  \ref{LC5}, the sequences
  $\left(\BU^n,(R'^n_p/\rn)_{p\geq1}\right)$
converge stably in law to $\left(\BU,(R'_p)_{p\geq1}\right)$ as $n\to\infty$.
\end{lem}

\nib Proof. \rm Due to Lemma \ref{LC5} and to the definition of
$R'_p$ and the fact that $\si$ is c\`adl\`ag, it is clearly enough to
prove that for any $p\geq1$ we have
\bee\label{C5}
w^n_p~:=~R'^n_p/\rn-\si_{S_-(n,p)}\al_-(n,p)-\si_{S_p}\al_+(n,p)~\toop~0.
\eee

We use the notations $\umu'$ and $(\f'_t)$ of the previous proof. We
deduce from (\ref{TT3}) that
\bee\label{TT33}
X(\ep)_t=X_0+\int_0^tb'(\ep)_sds+\int_0^t\si_sdW_s
+\de\star(\umu'-\unu')_t,
\eee
where $b'(\ep)t=b'_t-\int_E\de(t,x)dx$ and the above stochastic
integrals may be taken relative to both filtrations $(\f_t)$ and
$(\f'_t)$. In particular $X(\ep)$ satisfies (SH) for the filtration
$(\f'_t)$. Similar to $X'=X-X_0-X^c$, we write
$X'(\ep)=X(\ep)-X^c -X_0$. Then
$$w^n_p=\frac1{\rn}\left(\De^n_{I(n,p)}X'(\ep)
+\int_{S_-(n,p)}^{S_p}(\si_u-\si_{S(n,p)})dW_s
+\int_{S_p}^{S_+(n,p)}(\si_u-\si_{S_p})dW_s\right).$$

We may write (\ref{TR4}) for the process $X'(\ep)$ and with the
conditional expectations w.r.t.\
$\f'_{(i-1)\De_n}$ instead of $\f_{(i-1)\De_n}$. If we additionally
use the $\f'_0$-measurability of $I(n,p)$, and if we modify the
definition of $\chi'^n_i$ for $i=I(n,p)$ as to be the sum of the two
stochastic integrals in the previous display, we obtain by taking the
expectation :
\bean
\E(\phi_2(w^n_p))&\leq &K\rn+K
\E\Big(\frac1{\De_n}\int_{S_-(n,p)}^{S_+(n,p)}du
\int_{E^c\cap\{x:|\de(u,x)|\leq\De_n^{1/4}\}}\de(u,x)^2dx\Big)\\
&&+K\E\Big(\frac1{\De_n}\int_{S_-(n,p)}^{S_p}(\si_u-\si_{S(n,p)})^2du
.+\frac1{\De_n}\int_{S_p}^{S_+(n,p)}(\si_u-\si_{S_p})^2du\Big)
\eean
Since $|\de|\leq\ga$ and $\int\ga(x)^2dx<\infty$ and since $\si$ is
c\`adl\`ag bounded, we deduce from Lebesgue's theorem that the above goes
to $0$ as $n\to\infty$, hence $w^n_p\toop0$ and the result is proved. \qed

\subsection{Proof of Theorem \ref{T4}.}

We take $f\in\ea'_1\cap C^{0,\nu}$
and $\eta\in(0,\infty]$, with $\eta<\infty$ when $f$ is not bounded.
For any $\ep\in(0,\eta/2)$ we have the decomposition
$V^n(f)-\BH^n(f\psi_\eta)=\sum_{j=1}^3Z^n(\ep,j)$, where
$$Z^n(\ep,1)=V^n(f(1-\psi_\ep))-\BH^n(f\psi_\eta(1-\psi_\ep)),$$
$$Z^n(\ep,2)=Z^n(2)=\rn~\st\Big(|\be^n_i|-|m_1\si_{(i-1)\De_n}|\Big),$$
$$Z^n(\ep,3)_t=\st\Big(\ze^n_i(\ep)-\ec(\ze^n_i(\ep)\Big),\quad
\mbox{with}~~~\ze^n_i(\ep)=(f\psi_\ep)(\dd X)-|\rn~\be^n_i|.$$

First, Theorem \ref{T2} implies that $Z^n(\ep,1)\toSp 
\Si(f(1-\psi_\ep),\psi_\eta)$ because $f(1-\psi_\ep)\in\ea''_r\cap
C^{0,\nu}$ for all $r\in(1,2)$. Next, it
follows from the Lebesgue dominated convergence theorem for stochastic
integrals that
$\Si(f(1-\psi_\ep),\psi_\eta)\toucp\Si(f,\psi_\eta)$ as $\ep\to0$.
Lemma \ref{LCONV} implies that $Z^n(2)\tols
\sqrt{m_2-m_1^2}\int_0^t\si_udW'_u$. Hence,
due to the properties of the stable convergence in law, it
is enough to prove that for all $\rho>0$,
$$\lim_{\ep\to0}~\limsup_n~\PP\left(\sup_{s\leq t}|Z^n(\ep,3)_s|>\rho
\right)=0.$$
But $Z^n(\ep,3)$ is a locally square-integrable martingale w.r.t. the
filtration $(\f_{\De_n[t/\De_n]})_{t\geq0}$, whose predictable
quadratic variation $C^n$ satisfies $C^n_t\leq\st\ec(\ze^n_i(\ep)^2)$,
and for the above it is enough that
$$\lim_{\ep\to0}~\limsup_n~\PP\left(C^n_t>\rho\right)=0$$
for all $\rho>0$. This is obviously implied by (\ref{TT1}), and we
are done. \qed

\subsection{Proof of Theorem \ref{T5}.}

1) We first prove the result under the stronger assumption (SL-$s$).
Let $g$ be a $C^2_b$ and even function in general, or a $C^1$ even
function with $g'\in\ea$ when $X$ is continuous.  With the notation
(\ref{TR8'}) and (\ref{CD1}), we have
\bee\label{CD2}
\left.\begin{array}{l}
\frac1{\rn}\left(\De_nV'^n(g)_t-\int_0^t\rho_{\si_u}(g)\,du\right)=
\BU^n(g)_t+\frac1{\rn}~U^n(g)_t+M^n_t,\\[2mm]
\hskip5cm\mbox{where}~~~M^n_t=\st\Big(\ze^n_i-\ec(\ze^n_i)\Big),
\end{array}\right\}
\eee
and $\ze^n_i=\rn~\Big(g(\dd X/\rn)-g(\be^n_i)\Big)$. Since (SH)
holds, we can apply (\ref{TR7}) for $q=2$ to get
\bee\label{CC9}
\st\ec(|\ze^n_i|^2)~\toucp~0.
\eee
By Lenglart's inequality (as in Lemma \ref{LN1}) we deduce that
$M^n\toucp0$. Next, Lemma \ref{LCONV} implies that
$\BU^n(g)_t\tols \int_0^t\te_u\,d\BW_u$, where $\te_t=
\sqrt{\rho_{\si_t}(g^2)-\rho_{\si_t}(g)^2}$. Hence in view of
(\ref{CD2}) it remains to prove that $\frac1{\rn}~U^n(g)\toucp0$ when
$s\leq1$, and that $\De_n^{s/2-1}U^n(g)\toucp0$ when $s>1$\,: this is
given by Lemma \ref{LTT5}, case (i), and the result is proved.

2) Now we only assume (L-$s$). A
localization procedure, more sophisticated but similar to the one in
Lemma \ref{LN2}, is described in details in Section 3 of \cite{BGJPS}.
Namely, we find a sequence $(T_p)$ of stopping times increasing to
$+\infty$ and a sequence of processes $(X(p),\si(p))$, such that:

$\bullet$ If $X$ satisfies (L-$s$), then each $X(p)$ satisfies
(SL-$s$).

$\bullet$ We have $(X(p)_t,\si(p)_t)=(X_t,\si_t)$ for all $t<T_p$,
where $\si(p)$ denotes the process associated with $X(p)$ in
(\ref{AS3})-(\ref{AS4}).

Of course in \cite{BGJPS} this localization is done under the
additional assumption that $\de=0$ ($X$ is continuous), but the
presence of a non-vanishing $\de$ does not impair the argument.

On the one hand, Theorem \ref{T5} holds for each $X(p)$, that is:
$$\begin{array}{l}
s\leq1~~\Rightarrow~~\frac1{\rn}\left(\De_nV'^n(X(p);g)_t
-\int_0^t\rho_{\si(p)_u}(g)\,du\right)\tols
\int_0^t\te(p)_u\,d\BW_u\\[2mm]
s>1~~\Rightarrow~~\De_n^{1-s/2}\left(\De_nV'^n(X(p);g)_t
-\int_0^t\rho_{\si(p)_u}(g)\,du\right)\toucp 0,\end{array}$$
where $\te(p)_t=\sqrt{\rho_{\si(p)_t}(g^2)-\rho_{\si(p)_t}(g)^2}$. On
the other hand bothe the right and the left sides above, written for
$(X(p),\si(p)$ and also for $(X,\si)$ at time $t$, agree on the set
$\{t<T_p\}$. Since $T_p\to\infty$, we readily deduce Theorem \ref{T5}
for the initial process $X$.

\subsection{Proof of Theorem \ref{T6}.}

First, coming back to the localization procedure explained just above,
we have that if (L-$s$) and (H') hold for $X$, one can find the
processes $X(p)$ and $\si(p)$ such that (SL-$s$) {\em and}\ (SH')
hold. Then the same argument as in the previous proof shows that it is
enough to prove Theorem \ref{T6} under (SL-$s$) and (SH').

We take $f\in\ea_r$ for some $r\in(0,1]$. Instead of (\ref{CD2}) we have 
\bee\label{CD3}
\frac1{\rn}\left(\De_n^{1-r/2}V^n(f)_t-m_r\!\int_0^t\!\!c_u^{r/2}\,du\right)=
\BU^n(h_r)_t+\frac1{\rn}~U^n(h_r)_t+M^n_t+N^n_t,
\eee
where $M^n_t$ is like in (\ref{CD2}) with
$\ze^n_i=\rn~\Big(h_r(\dd X/\rn)-h_r(\be^n_i)\Big)$, and where
$N^n_t=\De^{1/2-r/2}V^n(f-h_r)$.

We have $|f-h_r|\leq k$ for some $K\in\ea'''_2\cap C^0$, hence
$|N^n|\leq \De_n^{1-r/2}V^n(k)$ and thus
\bee\label{CD4}
N^n~=~\left\{\begin{array}{ll} 
\tou(1)\quad&\mbox{if }~r<1\\ 
\tOu(1)\quad&\mbox{if }~r=1\end{array}\right.
\eee 
The other terms in (\ref{CD3}) are treated as in the proof of Theorem
\ref{T5}. We have $\BU^n(h_r)_t\tols$ 
$\sqrt{m_{2r}-m_r^2}\int_0^tc_u^{r/2}\,d\BW_u$ and $M^n\toucp0$ (we can
still apply (\ref{TR7}) with $q=2$ here to get (\ref{CC9})). Then in
view of (\ref{CD3}) and (\ref{CD4}) we readily deduce
Theorem \ref{T6} from Lemma \ref{LTT5}, case (iii).

Finally suppose that $X$ is continuous\,: we need to prove the result
without (SH'), when $r>1$. Since $\De_n^{1-r/2}V^n(h_r)=\De_nV'^n(h_r)$, 
it is a consequence of Theorem \ref{T5} when $f=h_r$, and
for $f\in\ea_r$ it
remains to prove that $\De_n^{1/2-r/2}V^n(f-h_r)\toucp0$. Note that
$|f-h_r|\leq Kh_q$ for some $K>0$ and some
$q>r-1$. Since $\E(|\dd X|^q)\leq K\De_n^{q/2}$ when $X$ is
continuous, we get
$$\E\left(\sup_{s\leq t}\left|\De_n^{1/2-r/2}V^n(f-h_r)_s\right|\right)
\leq K\De_n^{1/2-r/2}\st\E(|\dd X|^q)\leq Kt\De_n^{(1-r+q)/2},$$
hence the result.

\subsection{Proof of Theorem \ref{T6'}.}

Here again, as in the two previous proofs, it is enough to prove the
result under (SL-$s$). We have $\vpi\in(0,\frac12)$ and $\al>0$, and
we set $g_n=h_2\psi_{\al\De_n^\vpi}$ and
$\Bg_n=h_2\psi_{\al\De_n^\vpi/2}$. Since
$$\De_nV'^n(\Bg_n)~\leq~V''^n(\vpi,\al)~\leq~\De_nV'^n(g_n),$$
hence also
$|V''^n(\vpi,\al)-C|\leq|\De_nV'^n(g_n)-C|+|\De_nV'^n(\Bg_n)-C|$, it is 
clearly enough to prove the result for $\De_nV'^n(g_n)$ and $\De_nV'^n(\Bg_n)$
instead of $V''^n(\vpi,\al)$, and also that
\bee\label{CD6}
\rn~(V'^n(g_n)-V'^n(\Bg_n))~\toucp~0
\eee 
when $s\leq\frac{4\vpi-1}{2\vpi}$. 

Exactly as for (\ref{CD2}) we observe that, 
\bee\label{CD5}
\frac1{\rn}\left(\De_nV'^n(g_n)-C\right)=\BU^n(h_2)+\frac1{\rn}~U^n(g_n)+M^n,
\eee
where $M^n_t$ is like in (\ref{CD2}) with
$\ze^n_i=\rn~\Big(g_n(\dd X/\rn)-h_2(\be^n_i)\Big)$. Apply
(\ref{TT27}) to obtain $\st\ec((\ze^n_i)^2)~=~\tou(\De_n^{4\vpi-1-s\vpi})$,
hence by Lenglart inequality,
\bee\label{CD7}
M^n~=~\tou(\De_n^{2\vpi-1/2-s\vpi/2}).
\eee
Next, Lemma \ref{LCONV} implies $\BU^n(h_2)_t\tols\sqrt{2}\int_0^t
c_u\,d\BW_u$. Hence if we plug Lemma \ref{LTT5} for case (ii) and
(\ref{CD7}) into (\ref{CD5}) we obtain the desired results for
$\De_nV'^n(g_n)$ (observe that $2\vpi-1/2-s\vpi/2\geq0$ if
$s\leq(4\vpi-1)/2\vpi$, and otherwise
$2\vpi-1/2-s\vpi/2\geq2\vpi-1/2-s\vpi$), and of course the same holds
for $\De_nV'^n(\Bg_n)$. 

It remains to prove (\ref{CD6}) when $s\leq\frac{4\vpi-1}{2\vpi}$. We
have the same decomposition (\ref{CD5}) for $\Bg_n$, with some $\BM^n$
which also satisfies (\ref{CD7}), whereas the left side of (\ref{CD6})
is $M^n-\BM^n$, and this goes u.c.p.\ to $0$ by (\ref{CD7}), and the
proof is finished.

\subsection{Proof of Theorem \ref{T7}-(i).}

Our first task is to show that (\ref{C2}) makes sense, and for further purposes
we slightly extend the setting : Take any sequence
$(T_n)$ of stopping times whose graphs are pairwise disjoint, and such
that $\De X_t(\om)\neq0$ implies the existence of $n=n(\om,t)$ such
that $t=T_n(\om)$.

\begin{lem}\label{LC4} Under (H), and if $g\in\ea''_1$, the increasing process
$C(g)$ defined by (\ref{C2'}) is finite-valued, and the formula
(\ref{C2}) defines a semimartingale on the extended space
$\probt$, which is a locally square-integrable martingale as soon
as the process $C(g)$ above is locally integrable, in which case
\bee\label{C3'}
\E(Z(g)_T^2)=\E(C(g)_T)
\eee
for any $(\f_t)$-stopping time $T$. Moreover we have:

a) Conditionally on $\f$, the process $Z(g)$ is a square-integrable martingale
with independent increments and predictable bracket $C(g)$, relative
to the filtration $(\f\bigvee\Wf_t)$, and whose law is completely
characterized by the processes $X$ and $c$ and do not depend on the
particular sequence $(T_n)$ of stopping times.

b) If further $X$ and $\si$ have no common jumps, then conditionally
on $\f$ the process $Z(g)$ is a Gaussian martingale.
\end{lem}

In Theorem \ref{T7} we have a large degree of freedom for defining
$Z(f')$, its conditional law w.r.t.\ the $\si$-field $\f$ being the
only relevant property. This lemma shows that if we change the
sequence of stopping times $(T_n)$, subject of course to the property
of encompassing all jump times of $X$, then one changes $Z(f')$ but
{\em not}\, its conditional law. Note also that $g(0)=0$ above, so in
(\ref{C2}) the ``part'' of $T_n$ for which $\De X_{T_n}=0$ does not
come in into the sum, which is consistent with what precedes.
\vsq

\nib Proof. \rm Among several natural proofs, here is an ``elementary'' one.
Let $g\in\ea''_1$, and set $\al_n=g(X_{T_n})(c_{T_n-}+\frac12~\De c_{T_n})$.
We have $g^2\star\mu_t<\infty$ and $c$ is $\om$-wise
locally bounded, hence $C(g)_t=\sum_n\al_n1_{\{T_n\leq t\}}<\infty$
(for $\PP$-almost all $\om$  of course).

Fix $\om\in\Om$ such that $C(g)(\om)_t<\infty$ for all
$t<\infty$. Under $\PP'$, for all $n$ with $T_n(\om)\leq t$
the variables $A_n(\om):=\sqrt{\ka_n}~U_n~\si_{T_n-}(\om)+
\sqrt{1-\ka_n}~U'_n~\si_{T_n}(\om)$
are independent centered with variances $\al_n(\om)$. Then by a
standard
criterion for convergence of series of independent variables, the formula
$$Z(g)_t(\om,.)=\sum_{n=1}^\infty g(X_{T_n}(\om))\Big(
\sqrt{\ka_n}~U_n~\si_{T_n-}(\om)+
\sqrt{1-\ka_n}~U'_n~\si_{T_n}(\om)\Big)1_{\{T_n(\om)\leq t\}}$$
defines a process $(\om',t)\mapsto Z(g)(\om,\om')_t$ which obviously
is a  martingale with independent increments. Moreover its predictable
bracket is deterministic (that is, it does not depend on $\om'$) and is
$C(g)(\om)$, and it is purely discontinuous and jumps at times
$T_n(\om)$, and the law of the jump at $T_n(\om)<\infty$ is the law
of $g(X_{T_n}(\om))\Big(\sqrt{\ka_n}~U_n~\si_{T_n-}(\om)+
\sqrt{1-\ka_n}~U'_n~\si_{T_n}(\om)\Big)$, which only depends on the
processes $X$ and $c$ at point $\om$. If further $X$ and $c$ have no
common jumps, then the law of the jump at $T_n(\om)<\infty$ is the law
of $g(X_{T_n}(\om))~\si_{T_n}(\om)\Big(\sqrt{\ka_n}~U_n+
\sqrt{1-\ka_n}~U'_n\Big)$, which is $\n(0,g(X_{T_n}(\om))^2~c_{T_n})$.
This proves (a) and (b).

Next we consider the properties of $Z(g)$, considered now as a
process defined on $\probt$. Suppose first that
$\E(C(g)_{S_p})<\infty$ for some sequence $(S_p)$ of stopping times
increasing to $\infty$. Then for any $(\f_t)$-stopping time $T$ we
have
$$\E(Z(g)_T^2)=\int\PP(d\om)\int \PP'(d\om')Z(g)_{T(\om)}(\om,\om')^2
=\int\PP(d\om)C(g)_{T(\om)}(\om),$$
so (\ref{C3'}) holds and $Z(g)_{S_p\bigwedge t}^2$ is
$\WP$-integrable, and for $A\in\Wf_t$ and $s\geq 0$ we have
\bean
&&\E\left(1_A(Z(g)_{S_p\bigwedge(t+s)}-Z(g)_{S_p\bigwedge t})\right)\\
&&=\int\PP(d\om)\int \PP(d\om')
1_A(\om,\om')(Z(g)_{S_p(\om)\bigwedge(t+s)}(\om,\om')
-Z(g)_{S_p(\om)\bigwedge t}(\om,\om'))=0,
\eean
and thus $Z(g)$ is an $(\Wf_t)$-locally square-integrable
martingale. In the general case we set $A_n=\{\al_n\leq1\}$ and we let
$T'_n=T_n$ and $T''_n=\infty$ on $A_n$, and $T'_n=\infty$ and
$T''_n=T_n$ on $A_n^c$. These are stopping times, and we define
$Z'(g)_t$ and $Z''(g)_t$ by (\ref{C2}), with the sequences $(T'_n)$
and $(T''_n)$ respectively. The same analysis as above shows that
$Z'(g)$ is an $(\Wf_t)$-locally square-integrable martingale, whereas
$Z''(g)_t$ is a finite sum and thus as a process it has finite
variation. We deduce the semimartingale property of
$Z(g)=Z'(g)+Z''(g)$. \qed
\vsc

\begin{lem}\label{LC7} The claim (i) of Theorem \ref{T7} holds under
  (SK) and when $f$ is $C^1$ and vanishes on a neighborhood of $0$.
\end{lem}

\nib Proof. \rm We suppose that $f(x)=0$ if $|x|\leq2\ep$
for some $\ep>0$. We use the notation (\ref{C4})
associated with this particular $\ep$, so that $|\De X_s|\leq\ep$
identically if $s$ is not equal to one of the $S_p$'s. Since the
derivative $f'$ also vanishes on $[-2\ep,2\ep]$, we deduce that the process
$Z(f')$ has the same law, conditional on $\f$, than the following
process:
$$Z_t=\sum_{p:~S_p\leq t} f'(R_p)R'_p.$$
Hence the claim amounts to the stable convergence in law towards $Z'$,
for the sequence of processes $Z^n(f)/\rn$, where $Z^n(f)$ is given by
(\ref{KR0}).

Recall that $V(f)=f\star\mu$. In view of the properties of $f$ we
readily check that on the set $\Om_n(T,\ep)$ we have, for $t\leq T$ :
\bea
\frac1{\rn}~Z^n(f)_t&=&\frac1{\rn}
\sum_{p:~S_p\leq\De_n[t/\De_n]}\Big(f(R_p+R'^n_p)-f(R_p)\Big)\nonumber\\
&=&\sum_{p:~S_p\leq\De_n[t/\De_n]}f'(R_p+\widetilde{R}'^n_p)
\frac{R'^n_p}{\rn},\label{C11}
\eea
where $\widetilde{R}^n_p$ is between $R_p$ and $R_p+R'^n_p$.
Since $R^n_p\to 0$, hence $\widetilde{R}^n_p\to 0$ as well, and since
$f'$ is continuous and $\Om_n(T,\ep)\to\Om$, the result is a trivial
consequence of Lemma \ref{LC5}. \qed
\vsc

Now we can prove Theorem \ref{T7}-(i) under (SK). For each $\ep>0$,
we set $f_\ep=f\psi_\ep$, and Lemma \ref{LC7}
implies $Z^n(f-f_\ep)/\rn\tols Z(f'-f'_\ep)$. On the other
hand, $f'\in\ea''_1$ and thus $Z(f')$ exists and $C(f'_\ep)_t\to0$ pointwise
(Lebesgue's theorem, notation (\ref{C2'})) as $\ep\to0$, and
$C(f'_\ep)_t\leq K_t$, so (\ref{C3'}) and Doob's inequality yield
$Z(f'_\ep)\toucp0$ and thus $Z(f'-f'_\ep)\toucp Z(f')$ as
$\ep\to0$. Therefore it remains to prove the following :
\bee\label{C9}
\lim_{\ep\to0}~\limsup_n~ \PP\Big(\sup_{s\leq t}
\left|Z^n(f_\ep)_t/\rn~\right|>\eta\Big)=0,\quad\forall\eta>0,~\forall t>0.
\eee

Set
\bee\label{C19}
k_\ep(x,y)=f_\ep(x+y)-f_\eta(x)-f_\ep(y),\qquad
g_\ep(x,y)=k_\ep(x,y)-f'_\ep(x)y.
\eee
For $\ep$ small enough the function $f_\ep$ is $C^2$, and
$V(f)_\ep)=f_\ep\star\mu$ and (\ref{TT3}) holds, so It\^o's formula
yields that $Z^n(f_\ep)/\rn=A(n,\ep)^{(n)}+M(n,\ep)^{(n)}$, where
$M(n,\ep)$ is a locally square-integrable martingale, and with
\bee\label{C20}
A(n,\ep)_t=\int_0^ta(n,\ep)_u~du,\quad
A'(n,\ep)_t:=\langle M(n,\ep,M(n,\ep)\rangle=
\int_0^t a'(n,\ep)_u~du,
\eee
where
$$\left\{\begin{array}{l}
a(n,\ep)_t=\frac1{\rn}\Big(f_\ep'(X_t-X^{(n)}_t)b'_t
+\frac12~ f''_\eta(X_t-X^{(n)}_t)c_t
+\int g_\ep(X_t-X^{(n)}_t,\de(t,z))~dz\Big) \\[2mm]
a'(n,\ep)_t=\frac1{\De_n}\left(f'_\ep(X_t-X^{(n)}_t)^2c_t
+\int k_\ep(X_t-X^{(n)}_t,\de(t,z))^2~dz\right).  \end{array}\right.$$
In order to get (\ref{C9}), it is enough to prove the following, for
all $\eta>0$, $t>0$ :
\bee\label{C15}
\lim_{\ep\to0}~\limsup_n~\PP\Big(\sup_{s\leq t}(|A(n,\ep)_s|
+A'(n,\ep)_t>\eta\Big)~=~0.
\eee

Recall that $f(0)=f'(0)=0$ and $f''(x)=$ o$(|x|)$ as $x\to0$, so we have
\bee\label{C21}
j=0,1,2\quad\Rightarrow\quad|f_\ep^{j)}(x)|~\leq~\al_\ep~(|x|\bigwedge
\ep)^{3-j}
\eee
for
some $\al_\ep$ going to $0$ as $\ep\to0$, which implies
\bee\label{C22}
|k_\eta(x,y)|\leq K\al_\eta|x|~|y|,\qquad
|g_\eta(x,y)|\leq K\al_\eta|x|~y^2.
\eee
Then, in view of (SK), we deduce that $|a(n,\ep)_t|\leq
K\al_\ep|X_t-X^{(n)}_t|/\rn$ and $a'(n,\ep)_t|\leq
K\al_\ep|X_t-X^{(n)}_t|^2/\De_n$. Now, exactly
as for (\ref{TR2}), one readily checks that $\E(|X_{t+s}-X_t|^q)\leq
K_qs^{q/2}$ for all $q\in(0,2]$ and $s,t\geq0$, under (SH).
Applying this with $q=1$ and $q=2$, respectively, gives
$$ \E\left( v(A(n,\ep)_T\right)\leq KT\al_\ep,\qquad
\E\left(A'(n,\ep)_T\right)\leq KT\al_\ep^2,$$
and (\ref{C15}) immediately follows because $\al_\ep\to0$.
\vst

Finally it remains to prove the result under (K). This is
done using the same localization procedure than in Lemma \ref{LN2} or
in the proof of Theorems \ref{T5}, and we leave the
(easy) details to the reader.

\subsection{Proof of Theorem \ref{T7}-(ii).}

Before proceeding to the proof itself, we give two preliminary lemmas\,:
the first one is related to Lemma \ref{LT1}, the second one is a
simple application of It\^o's formula.

\begin{lem}\label{LC8} Under (SK) there exist increasing functions
$l_n$ on $(0,\infty)$ such that
\bee\label{C30}
\lim_{\eta\to0}~\limsup_n~l_n(\eta)~=~0,
\eee
and that for all $i,n\in\N$, $\ep,\eta>0$, we have
with $X(\ep)'=X(\ep)-X_0-X^c$)
\bee\label{C31}
t\leq\De_n\quad\Rightarrow\quad
\ec\left(|X(\ep)'_{(i-1)\De_n+t}-X(\ep)'_{(i-1)\De_n}|^2\bigwedge\eta^2
\right)~\leq~\De_nl_n(\eta).
\eee
\end{lem}

\nib Proof. \rm For any $\te>0$ we use the decomposition
$X(\ep)'=N(\te)+M(\te)+B (\te)$ given in the proof of Lemma \ref{LT1}, and
also the function $\ga_2(y)$ given in lemma \ref{LTT3}. Recall that
$$\left\{\begin{array}{l}
\PP^n_{i-1}(N(\te)_{(i-1)\De_n+t}-N(\te)_{(i-1)\De_n}\neq0)\leq K\te^{-2}t,\\
\ec((M(\te)_{(i-1)\De_n+t}-M(\te)_{(i-1)\De_n})\leq \ga_2(\te)t,\\
|B(\te)_{(i-1)\De_n+t}-B(\te)_{(i-1)\De_n}|\leq K\te^{-1}t
\end{array}\right.$$
and $K$ above does not depend on $\ep$.
The same argument than in Lemma \ref{LT1} shows that
$$\ec\left(|X(\ep)'_{(i-1)\De_n+t}-X(\ep)'_{(i-1)\De_n}|^2\bigwedge\eta^2
\right)~\leq~K\left(\frac{\eta^2\De_n}{\te^2}+\De_n\ga_2(\te)
+\frac{\De_n^2}{\te^2}\right),$$
as soon as $t\leq\De_n$. So we have (\ref{C31}) if we take
$l_n(\eta)=K\inf_{\te\in(0,1]}\Big(\eta^2\te^{-2}+\ga_2(\te)
+\De_n\te^{-2}\Big)$, which is obviously increasing in
$\eta$. Moreover we have (\ref{C30}), otherwise there would be an infinite
sequence $n_k$ and a number $a>0$ such that
$\eta^2\te^{-2}+\ga_2(\te)+\De_{n_k}\te^{-2}\geq a$ for all
$\te\in(0,1]$ and all $\eta>0$, and this contradicts the fact that
$\ga_2(\te)\to0$ as $\te\to0$. \qed

\begin{lem}\label{LC9} Under (SK) there is a constant $K_0$ such that,
for each $C^2$ function $g$ satisfying $g(0)=0$ and $|g'|\leq A$ and
$|g''|\leq A$, we have for all $i$, $n$ and all $\ep>0$ :
\bee\label{C32}
t\leq\De_n\quad\Rightarrow\quad\left\{\begin{array}{l}
\left|\ec\left(g(X(\ep)_{(i-1)\De_n+t}-X(\ep)_{(i-1)\De_n})\right)\right|
\leq K_0A\De_n,\\[2mm]
\ec\left(g(X(\ep)_{(i-1)\De_n+t}-X(\ep)_{(i-1)\De_n})^2\right)
\leq K_0(A+A^2)\De_n.
\end{array}\right.
\eee
If moreover (SL-$2$) holds we also have
\bee\label{C33}
t\leq\De_n\quad\Rightarrow\quad\left\{\begin{array}{l}
\left|\ec\left(c_{(i-1)\De_n+t}-c_{(i-1)\De_n}\right)\right|
\leq K\De_n,\\[2mm]
\ec\left(|c_{(i-1)\De_n+t}-c_{(i-1)\De_n}|^2\right)\leq K\De_n.
\end{array}\right.
\eee
\end{lem}

\nib Proof. \rm By (\ref{TT3}) and It\^o's formula, we have
\bean
g(X(\ep)_{(i-1)\De_n+t}-X(\ep)_{(i-1)\De_n})&\!\!=&\!\!\!
\int_{(i-1)\De_n}^{(i-1)\De_n+t}b(n,i,\ep)_u~du\!
+\!\!\int_{(i-1)\De_n}^{(i-1)\De_n+t}\si(n,i,\ep)_u~dW_u\\
&&\hskip.5cm
+\int_{(i-1)\De_n}^{(i-1)\De_n+t}\int_\R\de(n,i,\ep)(u,x)(\umu-\unu)(du,dx),
\eean
for suitable coefficients easy to compute and
which under (SK) satisfy
$$|b(n,i,\ep)_t|\leq KA,\quad |\si(n,i,\ep)_t|\leq KA,\quad
|\de(n,i,\ep)(t,x)|\leq KA\ga(x),$$
uniformly in all arguments (including $\om$...).
Then (\ref{C32}) follows in a classical way.

Under (SL-$2$) the process $\si_t$ satisfies (SK) (except that there
are two Brownian motions, but this makes no difference here), so
(\ref{C32}) applied with $g(x)=x^2$ yields (\ref{C33}). \qed
\vsc

Now we proceed to the proof of Theorem \ref{T7}-(ii). Upon using
the same localization procedure than in the proof of (i), we see that it is
enough to prove the result under (SL-$2$), which we assume thereon.
We suppose that $f\in\ea_2$, so for $\ep>0$ small enough the function
$f_\ep=f\psi_\ep$ is $C^\infty$ and coincides with $h_2\psi_\ep$. We
divide the proof into several steps. 
\vst

\noindent \it Step 1. \rm Fix $\ep>0$. We apply Lemma \ref{LC6} with
$d=1$ and $\BU^n=\BU^n(h_2)$ to obtain
$$\left(\BU^n_t,(R'^n_p/\rn)_{p\geq1}\right)
\tols\left(\sqrt{2}\int_0^tc_u~d\BW_u,(R'_p)_{p\geq1}\right).$$

On the one hand, the function $f-f_\ep$ satisfies (\ref{C11}), so the same
argument than in Lemma \ref{LC7} allows to deduce that
\bee
\left(\BU^n_t,\frac1{\rn}~Z^n(f-f_\ep)\right)
\tols\left(\sqrt{2}\int_0^tc_u~d\BW_u,Z(f'-f'_\ep)\right).\label{C23}
\eee

On the other hand, suppose for a while that $X$ is continuous. Then
both (\ref{TR7}) and (\ref{TT111}) hold for $g=h_2$, and so the proof
of Theorem \ref{T5} holds in this case as well, that is 
$$\frac1{\rn}~(\De_nV'^n(h_2)-C)-\BU^n=\frac1{\rn}~(V^n(h_2)-C)-
\BU^n\toucp0.$$
We also have 
\bean
\frac1{\rn}~\E(V^n((1-\psi_\ep)h_2)_t)&\leq&
\frac1{\rn}\st \E\left(|\dd X|^21_{\{|\dd X|>\ep\}}\right)\\
&\leq&\frac1{\ep^2\rn}\st \E\left(|\dd X|^4\right)~\leq~
Kt\rn 
\eean
because $\E(|\dd X|^4)\leq K\De_n^2$ by (\ref{TR2}) when $X$ is
continuous. Combining these two results yields
$\frac1{\rn}~(V^n(f_\ep)-C)-\BU^n(h_2)\toucp0$. Now, $X$ is
discontinuous, but applying what precedes to $X^c$ yields
$$\frac1{\rn}~\left(V^n(X^c;f_\ep)-C\right)-\rn~\BU^n~\toucp~0.$$
Hence (\ref{C23}) holds with $\BU^n$ substituted with
$\frac1{\rn}~\left(V^n(f_\ep;X^c)-C\right)$.
Since the stochastic integral process in the right side of (\ref{C23})
is continuous, we deduce that
$$\frac1{\rn}\left(V^n(X^c;f_\ep)_t-C_t+Z^n(f-f_\ep)_t\right)
\tols\sqrt{2}\int_0^tc_u~d\BW_u+Z(f'-f'_\ep).$$
Furthermore we have $Z(f'-f'_\ep)\toucp Z(f')$ as $\ep\to0$ (this is
like in the previous proof), whereas
$V^n(f)-V(f)^{(n)}=Z^n(f-f_\ep)+V^n(f_\ep)-C^{(n)}-f_\ep\star\mu$,
and also $V^n(f_\ep)_s=V^n(X(\ep);f_\ep)_s$ for all $s\leq t$ on the
set $\Om(t,\ep)$, which converges to $\Om$ as $\ep\to0$. Therefore, for
obtaining the result it remains to prove that
\bee\label{C25}
\lim_{\ep\to0}~\limsup_n~ \PP\Big(\sup_{t\leq T}
\left|Y^n(\ep)_t/\rn\right|>\eta\Big)=0,\quad\forall\eta>0,~\forall T>0.
\eee
where $Y^n(\ep)=V^n(X(\ep);f_\ep)-V^n(X^c;f_\ep)+C-C^{(n)}-
f_\ep\star\mu$.
\vst

\noindent \it Step 2. \rm Recall that for $\ep$ small enough the function
$f_\ep$ is $C^\infty$. Then It\^o's formula applied with (\ref{TT33}) yields
$Y^n(\ep)/\rn=A(n,\ep)^{(n)}+M(n,\ep)^{(n)}$, where $M(n,\ep)$ is a
locally square-integrable martingale, and we see that (\ref{C20}) holds with
$$\begin{array}{l}
a(n,\ep)_t=\frac1{\rn}\Big(\frac12~(f''_\ep(X(\ep)_t-X(\ep)^{(n)}_t)-
f''_\ep(X_t^c-(X^c)^{(n)}_t))c_t\\
\hskip2cm+f_\ep'(X(\ep)_t-X(\ep)^{(n)}_t)b'(\ep)_t+
\Bg_{\ep,t}(X(\ep)_t-X(\ep)^{(n)}_t)+(c_t-c_t^{(n)})\Big) \\[2mm]
a'(n,\ep)_t=\frac1{\De_n}\Big((f'_\ep(X(\ep)_t-X(\ep)^{(n)}_t)-
f'_\ep(X^c_t-(X^c)^{(n)}_t))^2c_t
+\Bk_{\ep,t}(X(\ep)_t-X(\ep)^{(n)}_t)\Big),  \end{array}$$
where we use the notation (\ref{C19})  and $E_\ep^c=\{x:\ga(x)\leq\ep\}$ and
$$\Bk_{\ep,t}(x)=\int_{E_\ep^c} k_\ep(x,\de(t,y))^2~dy,\qquad
\Bg_{\ep,t}(x)=\int_{E_\ep^c} g_\ep(x,\de(t,y))~dy.$$

Here again we are left to proving (\ref{C15}). This is more difficult
than for (i) of Theorem \ref{T7}, because (\ref{C21}) and (\ref{C22})
no longer hold. However $A'(n,\ep)$ is increasing, whereas
$|a(n,\ep)_t|\leq K$ because of (SL-$2$) and because the functions
$f'_\ep$, $f''_\ep$ and $\Bg_{\ep,t}$ are obviously bounded by a
constant not depending on $(\ep,t)$. Hence (\ref{C15}) will follow
if we prove that for all $\eta>0$, $t>0$ we have
\bee\label{C15'}
\lim_{\ep\to0}~\limsup_n~\PP\Big(\sup_{s\leq t}(|A(n,\ep)_s^{(n)}|
+A'(n,\ep)^{(n)}_t>\eta\Big)~=~0.
\eee
\vst

\noindent\it Step 3. \rm We will introduce below some decompositions for
$A(n,\ep)^{(n)}$ and $A'(n,\ep)^{(n)}$, namely
\bee\label{C26}
A(n,\ep)^{(n)}=\sum_{j=1}^6D^n(\ep,j),\quad
A'(n,\ep)^{(n)}=\sum_{j=7}^8D^n(\ep,j),
\eee
where $D^n(\ep,j)_t=\st\ze^n_i(\ep,j)$.
Then in order to get (\ref{C15'}) is it obviously enough to prove that
$\lim_{\ep\to0}~\limsup_n~\PP\Big(\sup_{s\leq t}
|D^n(\ep,j)_s|>\eta\Big)=0$ for each $j$. This property obviously holds if
\bee\label{C27}
\lim_{\ep\to0}~\limsup_n~ \E\Big(\st|\ze^n_i(\ep,j)|\Big)~=~0,
\eee
and it also holds if for all $\eta>0$ we have the following two
properties, as $n\to\infty$ :
\bee\label{C28}
\E\Big(\st|\ec(\ze^n_i(\ep,j))|\Big)~\to~0,\qquad
\E\Big(\st|\ze^n_i(\ep,j)|^2\Big)~\to~0.
\eee
\vst

\noindent\it Step 4. \rm Before deriving (\ref{C15'}) we state a
number of properties of the functions $f_\ep$, $\Bg_{\ep,t}$ and
$\Bk_{\ep,t}$ and their derivatives. These properties are elementary,
although sometimes tedious to derive, and they are based on the fact that
$f_\ep$ is $C^\infty$ for $\ep$ small enough, and $f_\ep(x)=x^2$
when $|x|\leq\ep$ and $f_\ep(x)=0$ when $|x|\geq2\ep$; we also use
(SK) for (\ref{C53}) below, where the notation $\ga_2(y)$ of
Lemma \ref{LTT3} is used. Here is the list of those properties :
\bee\label{C50}
|f_\ep^{(l)}(x)|\leq K_l\ep^{2-l}1_{\{|x|\leq2\ep\}},
\eee
\bee\label{C51}
|f'_\ep(x+y)-f'_\ep(x)|^2\leq K(x^4/\ep^2+y^2\bigwedge\ep^2),
\eee
\bee\label{C52}
|f''_\ep(x+y)-f''_\ep(x)|\leq K(x^2+y^2)/\ep^2,
\eee
\bee\label{C53}
|\Bg_{\ep,t}(x)|\leq K(x^2/\ep^2+|x|/\ep), \qquad
\Bk_{\eta,t}(x)\leq Kx^2\ga_2(\ep),
\eee
\bee\label{C54}
l=1,2\quad\Rightarrow\quad|\Bg^{(l)}_{\eta,t}(x)|\leq K\eta^{-l},
\eee
\bee\label{C55}
|\Bg_{\eta,t}(x)-\Bg_{\eta,s}(x)|\leq K|x|\int_{E_\ep^c}
|\de(t,z)-\de(s,z)|~\ga(z)~dz\leq K|x|\ga_2(\ep).
\eee
\vst

\noindent\it Step 5. \rm Now, recalling that the right limit $b'(\ep)_{t+}$
of $b'(\ep)$ exists, and with $\Bg_{\ep,t+}(x)=
\int g_\ep(x,\de_+(t,y))~dy$ (see the notation before (\ref{TT5}), we
set

$$\begin{array}{lll}
\ze^n_i(\ep,1)&=&\frac1{\rn}\itai(c_t-c_{(i-1)\De_n})~dt\\[2mm]
\ze^n_i(\ep,2)&=&\frac1{2\rn}\itai\Big(f''_\ep(X(\ep)_t-X(\ep)_{(i-1)\De_n})-
f''_\ep(X^c_t-X^c_{(i-1)\De_n})\Big)~c_t~dt\\[2mm]
\ze^n_i(\ep,3)&=&\frac1{\rn}~b'(\ep)_{(i-1)\De_n+}\itai
f'_\ep(X (\ep)_t-X (\ep)_{(i-1)\De_n})~dt\\[2mm]
\ze^n_i(\ep,4)&=&\frac1{\rn}~\itai
f'_\ep(X (\ep)_t-X(\ep)_{(i-1)\De_n})(b'(\ep)_t-
b'(\ep)_{(i-1)\De_n+})~dt\\[2mm]
\ze^n_i(\ep,5)&=&\frac1{\rn}~\itai
\Bg_{\ep,(i-1)\De_n+}(X (\ep)_t-X (\ep)_{(i-1)\De_n})~dt\\[2mm]
\ze^n_i(\ep,6)&=&\frac1{\rn}~\itai
\Big(\Bg_{\ep,t}-\Bg_{\ep,(i-1)\De_n+}\Big)(X(\ep)
_t-X (\ep)_{(i-1)\De_n})~dt\\[2mm]
\ze^n_i(\ep,7)&=& \frac1{\De_n}~\itai \Big(f'_\eta(X(\ep)_t
-X(\ep)_{(i-1)\De_n})-
f'_\ep(X^c_t-X^c_{(i-1)\De_n})\Big)^2~c_t~dt\\[2mm]
\ze^n_i(\ep,8)&=&\frac1{\De_n}~\itai
\Bk_{\ep,t}(X(\ep)_t-X(\ep)_{(i-1)\De_n})~dt
\end{array}$$
With these variables, it is easy to check that (\ref{C26})
holds. Hence it remains to prove that for each $j=1,\ldots,8$ we have
either (\ref{C27}) or (\ref{C28}). This is the aim of the following
lemma, which will end our proof.

\begin{lem}\label{LF1} We have (\ref{C27}) for $j=2,4,6,7,8$.
\end{lem}

\nib Proof. \rm Recalling $X(\ep)=X_0+X^c+X(\ep)'$, we deduce from
(\ref{C52}) that
$$|\ze^n_i(\ep,2)|\leq \frac{K}{\ep^2\rn}\itai
\left((X^c_t-X^c_{(i-1)\De_n})^2+(X(\ep)'_t-
X(\ep)'_{(i-1)\De_n})^2\right)~dt.$$
Applying (\ref{C32}) with $g(x)=x$ to the two processes $X^c$ and
$X(\ep)'$ readily gives $\ec(|\ze^n_i(\ep,2)|)\leq K\De_n^{3/2}/\ep^2$, and
(\ref{C27}) follows.

In a similar way, (\ref{C51}) gives
$$|\ze^n_i(\ep,7)|\leq \frac{K}{\De_n}\itai
\left(\ep^{-2}(X^c_t-X^c_{(i-1)\De_n})^4
+\Big((X(\ep)'_t-X(\ep)'_{(i-1)\De_n})^2\bigwedge\ep^2\Big)\right)~dt.$$
Applying the well known fact that
$\ec\Big((X^c_t-X^c_{(i-1)\De_n})^4\Big)\leq Kt^2$, and (\ref{C31}),
we deduce
$$\ec(|\ze^n_i(\ep,7)|)\leq K\left(\frac{\De_n^2}{\ep^2}
+\De_n l_n(\ep)\right).$$
Then we readily deduce (\ref{C27}) from (\ref{C30}).

Use (\ref{C55}) and (\ref{C53}), together with (\ref{C32}) again and
Cauchy-Schwarz for $j=6$, to get
$$ \ec(|\ze(\ep,6)|)+\ec(|\ze(\ep,8)|)\leq K\De_n\ga_2(\ep).$$
Since $\ga_2(\ep)\to0$ as $\ep\to0$,
we deduce (\ref{C27}) for $j=6,8$.

Finally consider the case $j=6$. We use (\ref{C50}) and (\ref{C32}) once
more, plus Cauchy-Schwarz, to get (with $b'^{(n)}_+$ being the process
associated with $(b'_{t+}$ by (\ref{I2}):
\bean
&&\E\left(\st|\ze^n_i(\ep,6)|\right)
\leq
\frac1{\rn}~\ec\left(\int_0^t|X(\ep)_s-X(\ep)^{(n)}_s|~
|b'(\ep)_s-b'(\ep)^{(n)}_{s+}|~ds \right)\\
&&\quad\leq
\frac1{\rn}\left( \ec\left(\int_0^t|X(\ep)_s-X(\ep)^{(n)}_s|^2
~ds\right)
\ec\left(\int_0^t
|b'(\ep)_s-b'(\ep)^{(n)}_{s+}|^2~ds\right)\right)^{1/2}\\
&&\quad\leq
\left(\ec\left(\int_0^t
|b'(\ep)_s-b'(\ep)^{(n)}_{s+}|^2~ds
\right)\right)^{1/2}
\eean
where the last inequality comes from (\ref{C32}). The last term above
goes to $0$ because $b'(\ep)_s-b'(\ep)^{(n)}_{s+}$ goes pointwise
to $0$ and is bounded: therefore we have (\ref{C27}) for $j=6$. \qed

\begin{lem}\label{LF2} We have (\ref{C28}) for $j=1,3,5$.
\end{lem}

\nib Proof. \rm Note that $\ze^n_i(\ep,1)=\ze^n_i(1)$ does not
depend on $\ep$. Then (\ref{C28}) for $j=1$ readily follows from
(\ref{C33}).

Next, use (\ref{C32}) for the function $f'_\ep$ and (\ref{C50}) and
the boundedness of $b'$ to obtain
$$\left|\ec(\ze^n_i(\ep,3))\right|\leq
\frac{K\De_n^{3/2}}{\ep},\qquad \ec(\ze^n_i(\ep,3)^2)\leq
\frac{K\De_n^2}{\ep^2},$$
and we readily deduce (\ref{C28}) for $j=3$. The same argument also
show (\ref{C28}) for $j=5$: we use (\ref{C54}), and (\ref{C32}) with
the function $\Bg_{\ep,(i-1)\De_n+}$ (this function is random, but
$\f_{(i-1)\De_n}$-measurable and with uniform bounds on its
derivatives, so (\ref{C32}) applies in this case). \qed

\subsection{Proof of Theorem \ref{T8}.}

Due to all what precedes, the proof is very easy\,: on the one hand,
Lemma \ref{LC6} is already multidimensional. On the other hand, the
way Theorems \ref{T5}, \ref{T6}, \ref{T6'} and \ref{T7} are deduced
from Lemma \ref{LC6} can be carried 
over separately for each component, in the multidimensional case.
Therefore Theorem \ref{T8} holds.

\end{document}